\documentclass[reqno,11pt]{amsart}
\usepackage{amssymb, amsmath,latexsym,amsfonts,amsbsy, amsthm}
\usepackage{color}

\let\pa=\partial
\let\al=\alpha

\let\f=\frac
\let\p=\psi

\let\Om=\Omega

\let\ve=\varepsilon
\let\pa=\partial

\def\na{\nabla}

\def\p{\partial}

\newcommand{\beq}{\begin{equation}}
\newcommand{\eeq}{\end{equation}}
\newcommand{\ben}{\begin{eqnarray}}
\newcommand{\een}{\end{eqnarray}}
\newcommand{\beno}{\begin{eqnarray*}}
\newcommand{\eeno}{\end{eqnarray*}}

\renewcommand{\theequation}{\thesection.\arabic{equation}}

\newtheorem{theorem}{Theorem}[section]

\newtheorem{proposition}[theorem]{Proposition}

\newtheorem{Theorem}{Theorem}[section]

\newtheorem{Remark}[Theorem]{Remark}

\newcommand{\ud}{\mathrm{d}}

\newcommand{\hh}{\mathbf{h}}

\newcommand{\mm}{\mathbf{m}}
\newcommand{\nn}{\mathbf{n}}

\newcommand{\vv}{\mathbf{v}}
\newcommand{\pp}{\mathbf{p}}

\newcommand{\xx}{\mathbf{x}}

\newcommand{\DD}{\mathbf{D}}
\newcommand{\FF}{\mathbf{F}}

\newcommand{\HH}{\mathbf{H}}
\newcommand{\II}{\mathbf{I}}

\newcommand{\MM}{\mathbf{M}}
\newcommand{\NN}{\mathbf{N}}
\newcommand{\QQ}{\mathbf{Q}}

\newcommand{\BS}{{\mathbb{S}^2}}
\newcommand{\BR}{{\mathbb{T}^3}}
\newcommand{\BT}{{\mathbb{R}^3}}

\newcommand{\BOm}{\mathbf{\Omega}}

\newcommand{\tr}{\mathrm{Tr}}

\begin{document}
\title[Nematic-Isotropic sharp interface]
{Dynamics of the Nematic-Isotropic sharp interface for the liquid crystal}

\author{Mingwen Fei}
\address{School of  Mathematics and Computer Sciences, Anhui Normal University, Wuhu, China}
\email{ahnufmwen@126.com}

\author{Wei Wang}
\address{BICMR, Peking University, Beijing 100871, China}
\email{wangw07@pku.edu.cn}

\author{Pingwen Zhang}
\address{School of  Mathematical Sciences, Peking University, Beijing 100871, China}
\email{pzhang@pku.edu.cn}

\author{Zhifei Zhang}
\address{School of  Mathematical Sciences, Peking University, Beijing 100871, China}
\email{zfzhang@math.pku.edu.cn}

\date{\today}
\maketitle

\renewcommand{\theequation}{\thesection.\arabic{equation}}
\setcounter{equation}{0}
\begin{abstract}
In this paper, we derive the sharp interface model of the nematic-isotropic phase transition from the Landau-de Gennes theory by using
the matched asymptotic expansion method. The model includes the evolution equation of the velocity and director field of the liquid crystal,
the sharp interface and Young-Laplace jump condition on the interface.
\end{abstract}

\numberwithin{equation}{section}

\indent

\section{Introduction}

Liquid crystals are a state of matter that have properties between those of a conventional liquid and those of a solid crystal.
They possess several phases, for example, the nematic phase at low temperature and isotropic phase (ordinary fluid) at high temperature.
In a nematic phase, the rod-like molecules have no positional order, but they self-align to have long-range directional order with their long axes roughly parallel. There is a first order phase transition to the nematic  liquid crystal phase at a clearing point $T_{NI}$.
There is also another temperature $T^*$. For $T<T^*$, the isotropic phase is unstable, and  the nucleation of the nematic  phases takes places through a phase-ordering kinetic process.
For $T^*<T<T_{NI}$, the isotropic phase is metastable.
A classical model which predicts such kind of  nematic-isotropic phase transition is the hard-rod model proposed by Onsager \cite{On}
(see \cite{LZZ} and references therein for rigorous results).
In this paper, we are concerned with the later region, in which stable nematic and isotropic phases can coexist.

There are two classical approaches to describe the nematic-isotropic interface. One approach uses the sharp interface model,
which involves solving the governing differential equations with matching boundary  conditions
at a moving interface separating the nematic phase and isotropic phase \cite{CFG,Fr}.
Another approach uses the phase-field model \cite{PSW, BBG}. The nematic-isotropic interface is modeled by a smooth transition region of finite width in the phase-field variable (or order parameter) $\QQ$.
The two phases correspond to regions in which the values of $\QQ$ are those pertaining to the isotropic and nematic phases.
In transition region, the order parameter varies continuously between its two equilibrium values.

In this paper, we choose the phase-field approach to study the dynamics of nematic-isotropic sharp interface
in the framework of  Landau-de Gennes  theory \cite{DG, MN}. In this theory,  the state of the nematic liquid crystals is described by the macroscopic Q-tensor order parameter,
which is a symmetric, traceless $3\times 3$ matrix. Physically, it can be interpreted as the second-order moment of
the orientational distribution function $f$, that is,
\beno
\QQ=\int_\BS(\mm\mm-\frac{1}{3}\II)f\ud\mm.
\eeno
When $\QQ=0$, the nematic liquid crystal is said to be isotropic. When $\QQ$ has two equal non-zero eigenvalues, it is said to be uniaxial and
$\QQ$ can be written as
\beno
\QQ=s\big(\nn\nn-\f13\II\big),\quad \nn\in \BS.
\eeno
When $\QQ$ has three distinct eigenvalues, it is said to be biaxial and $\QQ$ can be written as
\beno
\QQ=s\big(\nn\nn-\f13\II\big)+\lambda(\nn'\nn'-\frac13\II),\quad \nn,\,\nn'\in \BS,\quad \nn\cdot\nn'=0.
\eeno

The general Landau-de Gennes energy functional takes the form
\begin{align}
\nonumber\mathcal{F}(\QQ,\nabla\QQ)=&\int_{\Omega}\Big\{ \underbrace{\frac{a}2\tr\QQ^2
-\frac{b}{3}\tr\QQ^3+\frac{c}{4}(\mathrm{Tr}\QQ^2)^2}_{F_b:\text{bulk energy}}\\
&+\underbrace{\frac{1}{2}\Big(L_1|\nabla\QQ|^2+L_2Q_{ij,j}Q_{ik,k}
+L_3Q_{ij,k}Q_{ik,j}+L_4Q_{ij}Q_{kl,i}Q_{kl,j}\Big)}_{F_e:\text{elastic energy}} \Big\}\ud\xx.\label{eq:LG}
\end{align}
Here $\Om$ is a domain in $\BT$, $a, b, c$ are material-dependent and temperature-dependent nonnegative constants and $L_i(i=1,2,3,4)$ are material dependent elastic constants.
We refer to \cite{DG, MN} for more introduction.

Since the elastic constants $L_i(i=1,2,3,4)$ are typically very small compared with $a, b,c$,
we may introduce a small parameter $\varepsilon$ in (\ref{eq:LG}):
\begin{align}\label{eq:LG-ve}
\nonumber\mathcal{F}^\ve(\QQ,\nabla\QQ)=&\frac1{\ve^2}\int_\Omega\Big(\underbrace{\frac{a}2\tr\QQ^2
-\frac{ b}{3}\tr\QQ^3+\frac{c}{4}(\mathrm{Tr}\QQ^2)^2}_{F_b(\QQ)}\Big)\ud\xx
\\&+\int_\Omega\underbrace{\frac{1}{2}\Big(L_1|\nabla\QQ|^2+L_2Q_{ij,j}Q_{ik,k}
+L_3Q_{ij,k}Q_{ik,j}+L_4Q_{ij}Q_{kl,i}Q_{kl,j}\Big)}_{F_e(\QQ)} \ud\xx.
\end{align}
If $\QQ$ is uniaxial everywhere, i.e., $\QQ=s(\nn\nn-\frac13\II)$ with $\nn\in
\mathbb{S}^2$, then it is easy to see
\beno
F_b(\QQ)=\frac{s^{2}}{27}\big(9a-2bs+3cs^{2}\big),
\eeno
which has double well if and only if $b^2=27ac$. 
As the term $L_4Q_{ij}Q_{kl,i}Q_{kl,j}$ may cause the energy to be not bounded from below, we take $L_4=0$ in the sequel.

There are several dynamic Q-tensor models to describe the flow of the nematic liquid crystal,
which are either derived from the molecular kinetic theory for the rigid rods by various closure approximations
such as \cite{Feng, FLS}, or directly derived by variational method such as Beris-Edwards model \cite{BE} and Qian-Sheng's  model \cite{QS}.
In \cite{HLW}, we introduce a systematic schema to derive the continuum static and dynamic model from the molecular kinetic theory, where
we derive a dynamic Q-tensor model preserving the basic energy dissipation and the physical range of eigenvalues.
In this work, we will choose more popular Beris-Edwards model
\begin{align}
&\mathbf{v}^{\varepsilon}_t+\mathbf{v^{\varepsilon}\cdot\nabla v^{\varepsilon}}=-\mathbf{\nabla
p^{\varepsilon}}+\nabla\cdot(\sigma_{\varepsilon}^s+\sigma_{\varepsilon}^a+\sigma_{\varepsilon}^d),\label{BE-v}\\
&\mathbf{\mathbf{\nabla\cdot v^{\varepsilon}}}=0,\label{BE-d}\\
&\QQ^{\varepsilon}_t+\mathbf{v^{\varepsilon}}\cdot\nabla \QQ^{\varepsilon}+\QQ^{\varepsilon}\cdot\mathbf{\Omega^{\varepsilon}}
-\mathbf{\Omega^{\varepsilon}}\cdot\QQ^{\varepsilon}=\frac{1}{\Gamma}\HH^{\varepsilon}+\mathbf{S}_{\QQ^{\varepsilon}}(\DD^{\varepsilon}).\label{BE-Q}
\end{align}
Here $\mathbf{\mathbf{v^{\varepsilon}}}$ is the  velocity of the fluid, $\mathbf{\mathbf{p^{\varepsilon}}}$ is the pressure, $\Gamma$
is a collective rotational diffusion constant,
$\DD^{\varepsilon}=\frac{1}{2}(\nabla\mathbf{v^{\varepsilon}}+(\nabla\mathbf{v^{\varepsilon}})^{T})$,
$\mathbf{\Omega^{\varepsilon}}=\frac{1}{2}(\nabla\mathbf{v^{\varepsilon}}-(\nabla\mathbf{v^{\varepsilon}})^{T})$,
$\sigma_{\varepsilon}^s,\sigma_{\varepsilon}^a$ and
$\sigma_{\varepsilon}^d $ are symmetry viscous stress, anti-symmetry
viscous stress and distortion stress respectively defined by
\begin{align}
\sigma_{\varepsilon}^s=2\nu \mathbf{D}^{\varepsilon}-\mathbf{S}_{\QQ^{\varepsilon}}(\HH^{\varepsilon}),
\quad \sigma_{\varepsilon}^a=\QQ^{\varepsilon}\cdot\HH^{\varepsilon}-\HH^{\varepsilon}\cdot\QQ^{\varepsilon},
\quad (\sigma_{\varepsilon}^d)_{ij}=-\frac{\partial \mathcal{F}^\ve}{\partial Q_{kl,j}}Q_{kl,i}^\ve,
\end{align}
where $\nu>0$ is the viscous coefficient and $\HH^{\varepsilon}$ is the molecular field defined by
\begin{align}
\HH^\ve(\QQ^\ve)=&\frac{\delta \mathcal{F}^\ve}{\delta \QQ}
=-\f 1 {\ve^2}\frac{\partial {F}_b}{\partial\QQ}+\partial_i\Big(\frac{\partial {F}_e}{\partial\QQ_{,i}}\Big),\nonumber\\
=&-\f1 {\ve^2}f(\QQ)-\mathcal L \QQ,
\end{align}
where $f(\QQ)=a\QQ-b\QQ^2+c|\QQ|^2\QQ+\frac{b}{3}|\QQ|^2\II$ and
\beno
\big(\mathcal L\QQ\big)_{kl}=L_1\Delta Q_{kl}+\frac12(L_2+L_3)\big(Q_{km,ml}+Q_{lm,mk}-\frac23\delta_{kl}Q_{ij,ij}\big).
\eeno
And $\mathbf{S}_{\QQ^{\varepsilon}}(\MM)$ is defined by
\begin{align}
\mathbf{S}_{\QQ^{\varepsilon}}(\MM)=\xi\Big(\MM\cdot(\QQ^{\varepsilon}+\frac{1}{3}\II)+(\QQ^{\varepsilon}+\frac{1}{3}\II)\cdot\MM
-2(\QQ^{\varepsilon}+\frac{1}{3}\II)\MM:\QQ^{\varepsilon}\Big)
\end{align}
for symmetric and traceless matrix $\MM$, where $\xi$ is a constant
depending on the molecular details of a given liquid crystal.

To neglect the boundary effect, we consider the case of the domain $\Omega=\BR$. For the simplicity of notations,
we take $\nu=\Gamma=1$ in the sequel.

\section{Sharp interface model}

We will use the matched asymptotic expansion method  motivated by \cite{ABC} to study the behaviour of the solution $(\vv^\ve, \QQ^\ve)$ of the system (\ref{BE-v})--(\ref{BE-Q}) when $\ve $ is small.
The idea is to expand the solution in powers of $\ve$ away from the transition region(outer expansion) and inside the transition region(inner expansion).
By substituting these expansions into the equations and matching powers of $\ve$, one can determine the limit equation.
By matching the inner and outer expansions on the boundaries of the transition region, one can derive the jump condition on the sharp interface.

Now we present a sketch of our main results.
Assume that there exists a transition region of width $\ve$ separating two domains $\Om^\pm(t)$.
Let $\Gamma(t)$ be a smooth surface centered in the transition region and $\varphi$ be the signed distance to the sharp interface.

\begin{itemize}

\item In the region $\Om^\pm(t)$, the solution has the expansion in $\ve$:
\begin{align*}
&\vv^{\varepsilon}(t,x)=\vv^{(0)}_\pm(\tau,t,x)+\varepsilon\mathbf{v}^{(1)}_\pm(\tau, t,x)+\cdots,\\&
\mathbf{p}^{\varepsilon}(t,x)=\varepsilon^{-1}\mathbf{p}^{(-1)}_\pm(\tau, t,x)+\pp^{(0)}_\pm(\tau,t,x)+\cdots,\nonumber\\&
\QQ^{\varepsilon}(t,x)=\QQ^{(0)}_\pm(\tau,t,x)+\varepsilon\QQ^{(1)}_\pm(\tau,t,x)+\cdots.
\end{align*}
where $\tau=\frac t \varepsilon, \QQ^{(0)}=s_\pm(\nn\nn-\f13 \II)$ with $s_+=\frac{b+\sqrt{b^2-24ac}}{4c}\big(=\frac{b}{3c}=\frac{9a}{b}\big)$ and $s_-=0$.
The leading order term $\pp^{(-1)}$ of the pressure is harmonic in $\Om^\pm(t)$.
When $\tau\rightarrow +\infty$, $\big(\vv^{(0)}_+, \pp^{(0)}_+,\nn\big)$ satisfies
the Ericksen-Lesile system in $\Om^+(t)$:
\begin{align*}
&\pa_t\mathbf{v}^{(0)}_++\mathbf{v^{(0)}_+\cdot\nabla v^{(0)}_+}=-\mathbf{\nabla
p^{(0)}_+}+\nabla\cdot\big(\sigma^L+\sigma^E\big),\\
&\na\cdot\vv^{(0)}_+=0,\\
&\nn\times\big(-\Delta\nn+\NN-\DD^{(0)}_+\cdot\nn\big)=0.
\end{align*}
Here $\sigma^L$ is the Leslie stress and $\sigma^E$ is the Ericksen stress(see section 5.1).
While,$\big(\vv^{(0)}_-, \pp^{(0)}_-\big)$ satisfies the incompressible Navier-Stokes equations in $\Om^-(t)$:
\begin{align*}
&\pa_t\mathbf{v}^{(0)}_-+\mathbf{v^{(0)}_-\cdot\nabla v^{(0)}_-}=-\mathbf{\nabla
p^{(0)}_-}+\Delta \vv^{(0)}_-,\\
&\na\cdot\vv^{(0)}_-=0.
\end{align*}

\item In the transition region, the solution has the expansion in $\ve$:
\beno
&&\vv^\ve(t,x,z)=\widetilde{\vv}^{(0)}(t,x,z)+\ve\widetilde{\vv}^{(1)}(t,x,z)+\cdots,\\
&&\QQ^\ve(t,x,z)=\widetilde{\QQ}^{(0)}(t,x,z)+\ve\widetilde{\QQ}^{(1)}(t,x,z)+\cdots,
\eeno
with $z=\f {\varphi(t,x)} \ve$, where $\widetilde{\QQ}^{(0)}$ satisfies
\beno
\widetilde{\QQ}_{zz}^{^{(0)}}-f(\widetilde{\QQ}^{(0)})=0,\quad \widetilde{\QQ}^{^{(0)}}\rightarrow s_\pm(\nn\nn-\frac 13\II)\quad z\rightarrow\pm \infty,
\eeno
which has the uniaxial solution $s(z)(\nn\nn-\frac{1}{3}\II)$ with $s(z)$ satisfying
\beno
-s''+as-\frac{b}{3}s^2+\frac{2}{3}cs^3=0, \ s(-\infty)=0,\ s(+\infty)=s_+.
\eeno

\item  The sharp interface is described by the transported mean curvature flow:
\beno
\varphi_t-\Delta\varphi+\mathbf{v}^{(0)}\cdot\nabla\varphi=0.
\eeno

\item Jump conditions on $\Gamma(t)$:
\beno
&&\big[\vv^{(0)}\big]=0,\quad \nu\cdot\na \nn|_{\Gamma(t)}=0,\\
&&\big[\pp^{(-1)}\big]=-\frac{2}{3}\int_{-\infty}^{+\infty}\big|s'(z)\big|^2dzH,\quad \big[\pp^{(0)}\big]=\Big[\big\langle\sigma^L,\nu\otimes\nu\big\rangle\Big],
\eeno
where $\nu$ is the unit normal of $\Gamma(t)$, and $H$ is the mean curvature of $\Gamma(t)$.

\end{itemize}

\section{critical points of the bulk energy}

A matrix $\QQ_0$ is called a critical point of the bulk energy $F_b(\QQ)$ if $f(\QQ_0)=0$.
We have the following characterization for critical points \cite{WZZ}.

\begin{proposition}\label{prop:critical ponit}
$f(\QQ)=0$ if and only if
\begin{align*}
\QQ=s(\nn\nn-\frac13\II),
\end{align*}
for some $\nn\in \BS$, where $s=0$ or a solution of $2cs^2-bs+3a=0$, that is,
\beno
s_1=\frac{b+\sqrt{b^2-24ac}}{4c},\quad s_2=\frac{b-\sqrt{b^2-24ac}}{4c}.
\eeno
 Moreover, the critical point $\QQ_0=s_1\big(\nn\nn-\frac13\II\big)$ is stable.
\end{proposition}

Given a critical point $\QQ_0$, the linearized operator $f'(\QQ_0)$ of $f(\QQ)$ around $\QQ_0$ is given by
\begin{align}\label{equ:bulk-L}
f'(\QQ_0)\QQ&=a\QQ-b\big(\QQ_0\cdot\QQ+\QQ\cdot\QQ_0\big)+c|\QQ_0|^2\QQ+2(\QQ_0:\QQ)\big(c\QQ_0+\frac{b}{3}\II\big).
\end{align}
It is easy to compute that

\begin{align}\label{equ:bulk-S}
\big\langle f''(\QQ_0)\QQ_1, \QQ_2\big\rangle=\big\langle f''(\QQ_0)\QQ_2, \QQ_1\big\rangle=&-b\big(\QQ_1\cdot\QQ_2+\QQ_2\cdot\QQ_1\big)+2c \big(\QQ_0:\QQ_2\big)\QQ_1
\nonumber\\&+2c\big(\QQ_0:\QQ_1\big)\QQ_2+2\big(\QQ_1:\QQ_2\big)\big(c\QQ_0+\frac{b}{3}\II\big).
\end{align}

In view of Proposition 2.2 in \cite{WZZ}, we know that
\begin{proposition}\label{prop:kernel}
Let $\QQ_0=s(\nn\nn-\frac13\II)$ be a critical point with $s\neq 0$. Then the kernel space of the linearized  operator $f'(\QQ_0)$ is given by
\beno
\text{ker}f'(\QQ_0)=\Big\{\nn\nn^{\bot}+\nn^{\bot}\nn:\nn^{\bot}\in\mathbf{V_n}\Big\},
\eeno
 where $\mathbf{V_n}=\big\{\nn^{\bot}\in\mathbb{R}^3:\nn^{\bot}\cdot\nn=0\big\}$.
\end{proposition}

\section{Dynamics of sharp interface without hydrodynamics}

In this section, we consider the system without hydrodynamics. In such case, the system (\ref{BE-v})--(\ref{BE-Q}) is reduced to a gradient system of $\QQ^\ve$:
\ben\label{equ:LG-gradient}
\QQ^{\varepsilon}_t=-\frac{\delta \mathcal{F}^\ve}{\delta \QQ}=\f1 {\ve^2}f(\QQ^{\varepsilon})+\mathcal L \QQ^{\varepsilon}.
\een
Assume that there exists a transition region of width $\ve$ separating two domains $\Om^\pm(t)$.
Let $\Gamma(t)$ be a smooth surface centered in the transition region and $\varphi(t,x)$ be the signed distance to the sharp interface.

\subsection{Outer expansion}

We make a formal expansion for $\QQ^\ve$ in $\Om_\pm(t)$:
\begin{align}\label{equ:Q-exp}
\QQ^{\varepsilon}(t,x)=\QQ^{(0)}_\pm(t,x)+\varepsilon\QQ^{(1)}_\pm(t,x)+\varepsilon^2\QQ^{(2)}_\pm(t,x)+\cdots.
\end{align}
Then by Taylor expansion, we get
\begin{align}\label{equ:f-exp}
&f(\QQ^{\varepsilon})
=f(\QQ^{(0)}_\pm)+\varepsilon f'(\QQ^{(0)}_\pm)\QQ^{(1)}_\pm+\varepsilon^2\Big( f'(\QQ^{(0)}_\pm)\QQ^{(2)}_\pm
+\frac{1}{2}\big\langle f''(\QQ^{(0)}_\pm)\QQ^{(1)}_\pm,\QQ^{(1)}_\pm\big\rangle\Big)+O(\ve^3).
\end{align}

Plugging (\ref{equ:Q-exp}) and (\ref{equ:f-exp}) into (\ref{equ:LG-gradient}), then equating $\ve^k(k=-2,-1,0)$, we find that
\begin{align}
&f(\QQ^{(0)}_\pm)=0,\label{equ:Q-ve-2}\\
&f'(\QQ^{(0)}_\pm)\QQ^{(1)}_\pm=0,\label{equ:Q-ve-1}\\
&\pa_t\QQ_{\pm}^{(0)}=\mathcal L \QQ^{(0)}_\pm-f'(\QQ^{(0)}_\pm)\QQ^{(2)}_\pm-\frac{1}{2}\big\langle f''(\QQ^{(0)}_\pm)\QQ^{(1)}_\pm,\QQ^{(1)}_\pm\big\rangle.\label{equ:Q-ve0}
\end{align}

By Proposition \ref{prop:critical ponit}, the equation (\ref{equ:Q-ve-2}) ensures that $\QQ^{(0)}_\pm(t,x)=s_\pm\big(\nn(t,x)\nn(t,x)-\frac13\II\big)$ for
some $\nn(t,x)\in \BS$.
The equation (\ref{equ:Q-ve-1}) tells us that $\QQ^{(1)}_+ \in\text{ker}f'(\QQ^{(0)}_+)$. This means by Proposition \ref{prop:kernel} that
\begin{align}
\QQ^{(1)}_+(t,x)=\nn(t,x)\nn(t,x)^{\bot}+\nn(t,x)^{\bot}\nn(t,x),\quad \nn(t,x)^{\bot}\in \mathbf{V}_\nn,
\end{align}
from which, it follows that
\begin{align}\label{equ:Q-3.8}
\frac{1}{2}\big\langle f''(\QQ^{(0)}_+)\QQ^{(1)}_+,\QQ^{(1)}_+\big\rangle=&(2s_+-b)|\nn(t,x)^\perp|^2\nn(t,x)\nn(t,x)-b\nn(t,x)^{\bot}\nn(t,x)^{\bot}\nonumber\\
&+\frac{2}{3}(b-s_+)|\nn(t,x)^\perp|^2\II\perp\text{ker}f'(\QQ^{(0)}_+).
\end{align}

The solvability of (\ref{equ:Q-ve0}) to find $\QQ^{(2)}_+$ requires that
\begin{align}
\pa_t\QQ^{(0)}_+-\mathcal L \QQ^{(0)}_++\frac{1}{2}\big\langle f''(\QQ^{(0)}_+)\QQ^{(1)}_+,\QQ^{(1)}_+\big\rangle \perp\text{ker}f'(\QQ^{(0)}_+).
\end{align}
Thanks to (\ref{equ:Q-3.8}), it is enough to require that
\beno
\pa_t\QQ^{(0)}_+-\mathcal L \QQ^{(0)}_+\perp\text{ker}f'(\QQ^{(0)}_+).
\eeno
This means that
\beno
\big(\pa_t\QQ^{(0)}_+-\mathcal L \QQ^{(0)}_+\big):\big(\nn\nn^{\bot}+\nn^{\bot}\nn\big)=0
\eeno
for any $\nn^{\bot}\in  \mathbf{V}_\nn$. From the proof of Lemma 3.2 in \cite{WZZ}, we know that
\beno
&&\pa_t\QQ_+^{(0)}:\big(\nn\nn^{\bot}+\nn^{\bot}\nn\big)=2s_+\nn_t\cdot\nn^\bot,\\
&&-\mathcal L \QQ^{(0)}_+:\big(\nn\nn^{\bot}+\nn^{\bot}\nn\big)=\frac{1}{s_+}\hh\cdot\nn^\bot,
\eeno
where $\hh=-\frac{\delta E(\nn,\nabla\nn)}{\delta\nn}$ with $E(\nn,\nabla\nn)$ the Oseen-Frank energy defined by
\ben
E=\f {k_1} 2(\na\cdot\nn)^2+\f {k_2} 2\big(\nn{\cdot}(\na\times\nn)\big)^2
+\f {k_3} 2|\nn{\times}(\na\times \nn)|^2
+\frac{k_2+k_4}2\big(\textrm{tr}(\na\nn)^2-(\na\cdot\nn)^2\big),
\een
where the elastic constants $k_1, k_2, k_3, k_4$ are given by
\begin{align}\label{OF-LD-relation}
k_1=k_3=\big(2L_1+L_2+L_3\big)(s_+)^2,\quad k_2=2L_1(s_+)^2,\quad k_4=L_3(s_+)^2.
\end{align}

Hence, we conclude that $\nn$ satisfies
\ben\label{equ:Heatflow}
\big(2(s_+)^2\nn_t+\hh\big)\times \nn=0.
\een
In special case $L_1=1, L_2=L_3=0$(thus, $k_1=k_2=k_3=2(s_+)^2, k_4=0$), we have
$$
E=(s_+)^2|\na \nn|^2,\quad \hh=-2(s_+)^2\Delta\nn.
$$
In this case, the equation (\ref{equ:Heatflow}) is just the well-known harmonic heat flow.

\subsection{Inner expansion}

In the transition region, we make the following expansion for $\QQ^\ve$:
\begin{align}
\QQ^{\varepsilon}(t,x)=\widetilde{\QQ}^{(0)}(t,x,z)+\varepsilon\widetilde{\QQ}^{(1)}(t,x,z)+\varepsilon^{2}\widetilde{\QQ}^{(2)}(t,x,z)+\cdots,
\end{align}
where $z=\f {\varphi(t,x)} \ve$. Simple calculations give
\begin{align}
\QQ^{\varepsilon}_t=&\varepsilon^{-1}\varphi_t\widetilde{\QQ}_z^{(0)}+\varphi_t\widetilde{\QQ}_z^{(1)}+\widetilde{\QQ}_t^{(0)}+O(\ve),\nonumber\\
\mathcal L\QQ^{\varepsilon}=&\varepsilon^{-2}\mathcal A(\na\varphi, \widetilde{\QQ}^{(0)}_{zz})
+\varepsilon^{-1}\Big(\mathcal A(\na\varphi, \widetilde{\QQ}^{(1)}_{zz})+\mathcal B_1(\na\varphi, \na_x\widetilde{\QQ}_z^{(0)})
+\mathcal B_2(\na^2\varphi, \widetilde{\QQ}_z^{(0)})\Big)+O(1),\nonumber\\
f(\QQ^{\varepsilon})=&f(\widetilde{\QQ}^{(0)})+\varepsilon f'(\widetilde{\QQ}^{(0)})\widetilde{\QQ}^{(1)}
+\varepsilon^2\big(f'(\widetilde{\QQ}^{(0)})\widetilde{\QQ}^{(2)}+\frac{1}{2}\langle f''(\widetilde{\QQ}^{(0)})\widetilde{\QQ}^{(1)}, \widetilde{\QQ}^{(1)}\rangle\big)+O(\ve^3),\nonumber
\end{align}
where
\beno
&&\big(\mathcal A(\na\varphi, \QQ)\big)_{kl}=L_1{Q}_{kl}|\nabla\varphi|^2
+\f12(L_2+L_3)\big({Q}_{km}\p_m\varphi\p_l\varphi+{Q}_{lm}\p_m\varphi\p_k\varphi-\f23\delta_{kl}{Q}_{ij}\p_i\varphi\p_j\varphi\big),\\
&&\big(\mathcal B_1(\na\varphi, \na_x\QQ)\big)_{kl}=2L_1\pa_i{Q}_{kl}\pa_i\varphi+\f12(L_2+L_3)\big(\pa_m{Q}_{km}\pa_l\varphi
+\pa_l{Q}_{km}\pa_m\varphi+\pa_m{Q}_{lm}\pa_k\varphi\nonumber\\
&&\qquad\qquad+\pa_k{Q}_{lm}\pa_m\varphi-\f23\delta_{kl}\pa_i{Q}_{ij}\pa_j\varphi-\f23\delta_{kl}\pa_j{Q}_{ij}\pa_i\varphi\big),\nonumber\\
&&\big(\mathcal B_2(\na^2\varphi, \QQ)\big)_{kl}=L_1{Q}_{kl}\Delta\varphi+\f12(L_2+L_3)\big({Q}_{km}\pa_m\pa_l\varphi
+{Q}_{lm}\pa_m\pa_k\varphi-\f23\delta_{kl}{Q}_{ij}\pa_i\pa_j\varphi\big).
\eeno
Substituting the above expansion into (\ref{equ:LG-gradient}), then equating $\varepsilon^{-2},\varepsilon^{-1}$ terms, we obtain
\begin{align}
&-\mathcal A(\na\varphi, \widetilde{\QQ}^{(0)}_{zz})+f(\widetilde{\QQ}^{(0)})=0,\label{equ:Qin-2}\\
&-\mathcal A(\na\varphi, \widetilde{\QQ}^{(1)}_{zz})+f'(\widetilde{\QQ}^{(0)})\widetilde{\QQ}^{(1)}
=-\varphi_t\widetilde{\QQ}_z^{(0)}+\mathcal B_1(\na\varphi, \na_x\widetilde{\QQ}_z^{(0)})
+\mathcal B_2(\na^2\varphi, \widetilde{\QQ}_z^{(0)}).\label{equ:Qin-1}
\end{align}

Now let us derive the evolution equation of sharp interface $\varphi(t,x)$.
Assume that
\beno
\widetilde{\QQ}^{(0)}(t,x,z)\longrightarrow {\QQ}^{(0)}_\pm(t,x)\qquad \textrm{as  } z\rightarrow \pm \infty.
\eeno
Then by  integration by parts and (\ref{equ:Qin-2}), we find that
\ben
&&\int_{-\infty}^{+\infty}\Big(-\mathcal A(\na\varphi, \widetilde{\QQ}^{(1)}_{zz})+f'(\widetilde{\QQ}^{(0)})\widetilde{\QQ}^{(1)}\Big):\widetilde{\QQ}^{(0)}_zdz\nonumber\\
&&=\int_{-\infty}^{+\infty}\Big(-\mathcal A(\na\varphi, \widetilde{\QQ}^{(0)}_{zz})+f'(\widetilde{\QQ}^{(0)})\widetilde{\QQ}^{(1)}\Big):\widetilde{\QQ}^{(1)}_zdz=0.\label{equ:Q-3.16}
\een
Here we used the fact that $\widetilde{\QQ}^{(0)}, \widetilde{\QQ}^{(1)}$ are traceless and $\FF_b({\QQ}^{(0)}_\pm)=0$.
On the other hand, we have
\beno
&&\int_{-\infty}^{+\infty}\Big(-\varphi_t\widetilde{\QQ}_z^{(0)}+\mathcal B_1(\na\varphi, \na_x\widetilde{\QQ}_z^{(0)})
+\mathcal B_2(\na^2\varphi, \widetilde{\QQ}_z^{(0)})\Big):\widetilde{\QQ}_z^{(0)}dz\\
&&=-c\varphi_t+\na\cdot\big(A\na\varphi\big),
\eeno
where
\beno
&&c(t,x)=\int_{-\infty}^{+\infty}|\widetilde{\QQ}_z^{(0)}(t,x,z)|^2dz,\\
&&A_{kl}(t,x)=L_1c(t,x)\delta_{kl}+(L_2+L_3)\int_{-\infty}^{+\infty}\widetilde{\QQ}_{km,z}^{(0)}(t,x,z)\widetilde{\QQ}_{ml,z}^{(0)}(t,x,z)dz.
\eeno
This combined with (\ref{equ:Qin-1}) and (\ref{equ:Q-3.16}) gives
\ben
c\varphi_t-\na\cdot\big(A\na\varphi\big)=0.
\een

In special case $L_1=1, L_2=L_3=0$, we have
\ben\label{equ:Q-3.18}
c\varphi_t-\na\cdot\big(c\na\varphi\big)=0,
\een
and $\widetilde{\QQ}^{(0)}$ satisfies
\ben\label{equ:Q-3.19}
-\widetilde{\QQ}^{(0)}_{zz}+f(\widetilde{\QQ}^{(0)})=0
\een
together with the boundary conditions
\beno
&&\widetilde{\QQ}^{(0)}(t,x,z)\longrightarrow s_+\big(\nn(t,x)\nn(t,x)-\f13\II\big) \quad \textrm{as}\quad z\rightarrow +\infty,\\
&&\widetilde{\QQ}^{(0)}(t,x,z)\longrightarrow 0 \quad \textrm{as}\quad z\rightarrow -\infty.
\eeno
If the solution of (\ref{equ:Q-3.19}) is uniaxial with the form $s(z)\big(\nn(t,x)\nn(t,x)-\f13\II\big)$, then $s$ should satisfy
\begin{equation}\label{equ:s}
-s''+as-\frac{b}{3}s^2+\frac{2}{3}cs^3=0, \ s(-\infty)=0,\ s(+\infty)=s_+.
\end{equation}
Therefore, the function $c(t,x)$ is independent of $(t,x)$, and the equation (\ref{equ:Q-3.18}) is reduced to
\beno
\varphi_t-\Delta \varphi=0,
\eeno
which is the well-known mean curvature flow. Thus, the equation (\ref{equ:Qin-1}) is reduced to
\begin{align}
-\widetilde{\QQ}^{(1)}_{zz}+f'(\widetilde{\QQ}^{(0)})\widetilde{\QQ}^{(1)}=2\nabla\varphi\cdot\nabla\widetilde{\QQ}^{(0)}_z.\nonumber
\end{align}
Note that it holds for any $\widetilde{\mathbf{T}}=s(z)(\nn\nn'+\nn'\nn)$ with $\nn'\bot\nn$,
\begin{align}
-\widetilde{\mathbf{T}}_{zz}+\langle f'(\widetilde{\QQ}^{(0)}),\widetilde{\mathbf{T}}\rangle=0,\nonumber
\end{align}
which implies that
\begin{align}
\int_{-\infty}^{+\infty}\nabla\varphi\cdot\nabla\partial_z\big[s(z)(\nn\nn-\frac13\II)\big]:s(z)(\nn\nn'+\nn'\nn) dz=0.\nonumber
\end{align}
While, this is equivalent to
\begin{align}
\nabla\varphi\cdot\nabla(\nn\nn-\frac13\II):(\nn\nn'+\nn'\nn)=0,\nonumber
\end{align}
that is,
\begin{align}
\big(\na\varphi\cdot\nabla\nn\big)\cdot\nn'=0,\quad \text{ for all}\quad \nn'\bot\nn.\nonumber
\end{align}
On the other hand, we have
\begin{align}
\big(\na\varphi\cdot\nabla\nn\big)\cdot\nn=0.\nonumber
\end{align}
Especially, this means that $\nn$ should satisfy the Neumann condition on the sharp interface
\begin{align}
\nu\cdot\nabla\nn=0\quad \textrm{on}\quad \Gamma(t).\label{equ:n-neumann}
\end{align}
Here $\nu$ is the unit normal of $\Gamma$.

In summary, we derive the sharp interface model without hydrodynamics from the gradient system (\ref{equ:LG-gradient}):
\begin{align*}
&\nn_t-\Delta\nn=|\na \nn|^2\nn \quad \textrm{\textrm{in}} \quad \Omega^+(t),\\
&\nu\cdot\nabla\nn=0\quad \textrm{on}\quad \Gamma(t).
\end{align*}
While, the sharp interface $\Gamma(t)$ is determined by the mean curvature flow
\beno
\varphi_t-\Delta \varphi=0.
\eeno

\subsection{Asymptotic analysis of Landau-de Gennes energy}

Assume that the leading order term of $\QQ^\ve$ is is uniaxial with the form
\beno
\QQ^\ve(t,x)\sim s\big(\f {\varphi(t,x)} \ve\big)\big(\nn(t,x)\nn(t,x)-\f13\II\big)\triangleq \QQ_0(t,x).
\eeno
This is a reasonable assumption  at least in some special cases(for example, $L_2=L_3=0$ or $L_2<0, L_3=0$)
by the outer-inner asymptotic analysis and \cite{PWZZ}.

Without loss of generality we consider in the case of $a=\f13, b=3, c=1$ and then $s_+=1$. And we assume  $L_3=L_4=0$. Then  Landau-de Gennes energy takes
\beno
\mathcal{F}^\ve(\QQ,\nabla\QQ)=\frac1{\ve^2}\int_\BR\frac{a}2\tr\QQ^2
-\frac b3\tr\QQ^3+\frac{c}{4}(\mathrm{Tr}\QQ^2)^2\ud\xx+\frac 12\int_\BR L_1|\nabla\QQ|^2+L_2Q_{ij,j}Q_{ik,k}\ud\xx.
\eeno
Let $S_\ve(x)=s\big(\f {\varphi(t,x)} \ve\big)$. Putting $\QQ_0$ into $\mathcal{F}^\ve(\cdot,\nabla \cdot)$, we deduce that
\begin{align*}
\ve\mathcal{F}^\ve(\QQ_0,\nabla\QQ_0)
=&\frac{c}{9\ve}\int_\BR S^2_\ve(s_+-S_\ve)^2\ud\xx+\frac{\ve L_1}{2}\int_\BR \Big(\frac23|\nabla S_\ve|^2+S^2_\ve|\nabla\nn|^2\Big)\ud\xx\\
&+\frac{\ve L_2}{2}\int_\BR\Big(\frac19|\nabla S_\ve|^2+\frac13(\nn\cdot\nabla S_\ve)^2+S^2(|\nn\cdot\nabla\nn|^2+(\nabla\cdot\nn)^2)\\
&+\frac23S_\ve(2(\nn\cdot\nabla S_\ve)(\nabla\cdot\nn)-\nabla S_\ve\cdot(\nn\cdot\nabla\nn))\Big)\ud\xx\\
=&\frac{1}{9}\int_\BR \Big(\frac{cS^2_\ve(s_+-S_\ve)^2}{\ve}+\frac\ve 2(6L_1+L_2)|\nabla S_\ve|^2+\frac{3\ve L_2}{2}(\nn\cdot\nabla S_\ve)^2\Big)\ud\xx\\
&+\frac{\ve}{2}\int_\BR \Big(L_1S^2_\ve|\nabla\nn|^2+L_2S^2_\ve(|\nn\cdot\nabla\nn|^2+(\nabla\cdot\nn)^2)\Big)\ud\xx\\
&+\frac{\ve L_2}{2}\int_\BR\Big(\frac23S_\ve\big[2(\nn\cdot\nabla S_\ve)(\nabla\cdot\nn)-\nabla S_\ve\cdot(\nn\cdot\nabla\nn)\big]\Big)\ud\xx\\
\triangleq& A+B+C.
\end{align*}
Direct calculations lead to
\begin{align*}
A=&\frac1{9}\int_\Gamma\int_{-\infty}^{+\infty} \Big(c{s^2(s_+-s)^2}+\frac1 2(6L_1+L_2)s_z^2
+\frac{3 L_2}{2}(\nn\cdot\nu)^2 s_z^2\Big)\ud z\ud\sigma\\
=&\int_\Gamma(\alpha+\beta (\nn\cdot\nu)^2)\ud\sigma,\\
B=&\frac{\ve}{2}\int_{\Omega^+} \Big(L_1|\nabla\nn|^2+L_2(|\nn\cdot\nabla\nn|^2+(\nabla\cdot\nn)^2)\Big)\ud\xx+o(\ve),\\
C=&\frac{\ve L_2}{2}\int_\Gamma\int_{-\infty}^{+\infty}\Big(\frac23s\big[2(\nn\cdot\nu)s_z(\nabla\cdot\nn)-s_z\nu\cdot(\nn\cdot\nabla\nn)\big]\Big)\ud z\ud \sigma\\
=&\frac{\ve L_2}{6}\int_\Gamma\Big(2(\nn\cdot\nu)(\nabla\cdot\nn)-\nu\cdot(\nn\cdot\nabla\nn)\Big)\ud \sigma,
\end{align*}
where $\nu$ is the unit normal of $\Gamma$ and
\begin{align*}
&\alpha=\frac1{9}\int_{-\infty}^{+\infty} \Big(c{s^2(s_+-s)^2}+\frac1 2(6L_1+L_2)s_z^2\Big)\ud z,\\
&\beta=\frac{3 L_2}{18} \int_{-\infty}^{+\infty}s_z^2\ud z.
\end{align*}
This gives the following asymptotic of $\mathcal{F}^\ve(\QQ,\nabla\QQ)$ as $\ve\rightarrow 0$,
\begin{align}
\ve\mathcal{F}^\ve(\QQ_0,\nabla\QQ_0)=&\int_\Gamma(\alpha+\beta (\nn\cdot\nu)^2)\ud\sigma\nonumber\\
&+\ve\bigg\{\f12\int_{\Omega^+} \Big(L_1|\nabla\nn|^2+L_2(|\nn\cdot\nabla\nn|^2+(\nabla\cdot\nn)^2)\Big)\ud\xx\nonumber\\
&\qquad+\frac{L_2}{6}\int_\Gamma\Big(2(\nn\cdot\nu)(\nabla\cdot\nn)-\nu\cdot(\nn\cdot\nabla\nn)\Big)\ud \sigma\bigg\}+o(\ve),\label{equ:energy-exp}
\end{align}
which is consistent with the total free energy introduced in \cite{Rey}.
The second part on the right hand side corresponds to the well-known Oseen-Frank energy. Let us give some explanations for the other parts
from the point of view of energy minimization.

\begin{Remark}
We view $\al$ as a functional of $s$, then minimize $\al(s)$ in the class $s\in C^2, s(+\infty)=1, s(-\infty)=0$. That is,
\beno
\min_{s\in C^2, s(+\infty)=1, s(-\infty)=0}\frac1{9}\int_{-\infty}^{+\infty} \Big(c{s^2(s_+-s)^2}+\frac1 2(6L_1+L_2)s_z^2\Big)\ud z.
\eeno
If $s(z)$ is a minimizer, it should satisfy the following Euler-Lagrangian equation
\begin{equation}
-(6L_1+L_2)s''+2c(s_+^2s-3s_+s^2+2s^3)=0.\nonumber
\end{equation}
This is consistent with (\ref{equ:s}) in the case of $L_1=1, L_2=0$.

\end{Remark}

\begin{Remark}
For the case of $L_2<0$(hence, $\beta<0$), we know from (\ref{equ:energy-exp}) that $\nn$ should take the normal $\nu$ on $\Gamma$
in order that the energy is small as soon as possible.
In this case, we have by  $\nabla\cdot\nn=\nabla_\Gamma\cdot\nn+\nu\cdot\partial_\nu\nn$ that
 \begin{align*}
\int_\Gamma\Big(2(\nabla\cdot\nn)-\nu\cdot(\partial_\nu\nn)\Big)\ud \sigma=&\int_\Gamma\Big(2(\nabla_\Gamma\cdot\nn)+\nu\cdot(\partial_\nu\nn)\Big)\ud \sigma\\
=&4\int_\Gamma H\ud \sigma,
\end{align*}
where $H$ is the mean curvature of the interface $\Gamma$. If $L_2>0$, it seems reasonable to conjecture from (\ref{equ:energy-exp}) that $\nn$ should be tangent to $\Gamma(t)$.
However, the above analysis does not work, since the uniaxial solution is unstable in the case of $L_2>0$ by numerical analysis \cite{KB} and \cite{PWZZ}.
\end{Remark}

\section{Dynamics of sharp interface with hydrodynamics}

In this section, we will consider the case of $L_1=1, L_2=L_3=0$ in order to simplify the analysis. So, the molecular field $\HH^\ve=-\f1 {\ve^2}f(\QQ)-\Delta\QQ$.
Again, assume that there exists a transition region of width $\ve$ separating two domains $\Om^\pm(t)$.
Let $\Gamma(t)$ be a smooth surface centered in the transition region and $\varphi(t,x)$ be the signed distance to the sharp interface.

\subsection{Outer expansion}
In order to match strong singularity of the distortion stress $\sigma_\ve^d$, we need to introduce a fast time scale $\tau=\f t\ve$.
We make a formal expansion for $(\vv^\ve, \pp^\ve, \QQ^\ve)$ in $\Om_\pm(t)$:
\begin{align*}
&\vv^{\varepsilon}(t,x)=\vv^{(0)}_\pm(\tau,t,x)+\varepsilon\mathbf{v}^{(1)}_\pm(\tau, t,x)+\varepsilon^{2}\mathbf{v}^{(2)}_\pm(\tau,t,x)+\cdots,\\&
\mathbf{p}^{\varepsilon}(t,x)=\varepsilon^{-2}\mathbf{p}^{(-2)}_\pm(\tau,t,x)+ \varepsilon^{-1}\mathbf{p}^{(-1)}_\pm(\tau, t,x)+\pp^{(0)}_\pm(\tau,t,x)+\cdots,\nonumber\\&
\QQ^{\varepsilon}(t,x)=\QQ^{(0)}_\pm(\tau,t,x)+\varepsilon\QQ^{(1)}_\pm(\tau,t,x)+\varepsilon^2\QQ^{(2)}_\pm(\tau,t,x)+\cdots.
\end{align*}
Simple calculations lead to
\begin{align*}
&\HH^{\varepsilon}=-\varepsilon^{-2}f(\QQ^{(0)}_\pm)-\varepsilon^{-1}f'(\QQ^{(0)}_\pm)\QQ^{(1)}_\pm+\Big(\Delta\QQ^{(0)}_\pm-f'(\QQ^{(0)}_\pm)\QQ^{(2)}_\pm\nonumber\\
&\qquad-\frac{1}{2}\langle f''(\QQ^{(0)}_\pm)\QQ^{(1)}_\pm,\QQ^{(1)}_\pm\rangle\Big)+\cdots\triangleq \ve^{-2}\HH^{(-2)}_\pm+\ve^{-1}\HH^{(-1)}_\pm+\cdots,\\
&\mathbf{S}_{\QQ^{\varepsilon}}(\HH^{\varepsilon})=\varepsilon^{-2}\mathbf{S}_{\QQ^{(0)}_\pm}(\HH^{(-2)}_\pm)
+\sum\limits_{m=-1}^{0}\varepsilon^m\xi\bigg(\sum\limits_{k=-2}^{m}\HH^{(k)}_\pm\cdot\QQ^{(m-k)}_\pm+\sum\limits_{k=-2}^{m}\QQ^{(m-k)}_\pm\cdot\HH^{(k)}_\pm
+\frac{2}{3}\HH^{(m)}_\pm\nonumber\\
&\qquad-2\sum\limits_{k=0}^{m+2}\QQ^{(k)}_\pm\sum\limits_{l=-2}^{m-k}\HH^{(l)}_\pm:\QQ^{(m-k-l)}_\pm
-\f23\II\sum\limits_{k=-2}^{m}\HH^{(k)}_\pm:\QQ^{(m-k)}_\pm\bigg)+\cdots\nonumber\\
&\qquad\qquad\triangleq \varepsilon^{-2}\mathbf{S}_{\HH,\pm}^{(-2)}+\varepsilon^{-1}\mathbf{S}_{\HH,\pm}^{(-1)}+\cdots,\\
\end{align*}
and
\begin{align*}
&\sigma_{\varepsilon}^s=-\varepsilon^{-2}\mathbf{S}_{\HH,\pm}^{(-2)}-\varepsilon^{-1}\mathbf{S}_{\HH,\pm}^{(-1)}
+\big(\DD^{(0)}_\pm-\mathbf{S}_{\HH,\pm}^{(0)}\big)+\cdots \triangleq \varepsilon^{-2}\sigma_{(-2),\pm}^s+\varepsilon^{-1}\sigma_{(-1),\pm}^s+\cdots,\\
&\sigma_{\varepsilon}^a=\sum\limits_{m=-2}^{0}\varepsilon^m\bigg(\sum\limits_{k=-2}^{m}\QQ^{(m-k)}_\pm\cdot\HH^{(k)}_\pm
-\sum\limits_{k=-2}^{m}\HH^{(k)}_\pm\cdot\QQ^{(m-k)}_\pm\bigg)+\cdots\\
&\quad\triangleq \varepsilon^{-2}\sigma_{(-2),\pm}^a+\varepsilon^{-1}\sigma_{(-1),\pm}^a+\cdots,\\&
\sigma_{\varepsilon}^d=\sigma^d(\QQ^{(0)}_\pm,\QQ^{(0)}_\pm)+\cdots\triangleq \sigma_{(0),\pm}^d+\cdots.
\end{align*}

Plugging these expansions into (\ref{BE-v})-(\ref{BE-Q}), then equating the $\varepsilon^k(k=-2,-1,0)$ terms in (\ref{BE-v}), we get
\begin{align}
&\nabla\mathbf{ p^{(-2)}_\pm}=\nabla\cdot\sigma_{(-2),\pm}^s+\nabla\cdot\sigma_{(-2),\pm}^a,\label{equ:p-2-4}\\
&\pa_\tau\mathbf{v}^{(0)}_\pm=-\nabla\mathbf{p^{(-1)}_\pm}+\nabla\cdot\sigma_{(-1),\pm}^s+\nabla\cdot\sigma_{(-1),\pm}^a,\label{equ:p-1-4}\\
&\pa_\tau\mathbf{v}^{(1)}_\pm+\pa_t\mathbf{v}^{(0)}_\pm+\mathbf{v^{(0)}_\pm\cdot\nabla v^{(0)}}_\pm=-\mathbf{\nabla
p^{(0)}_\pm}+\nabla\cdot\big(\sigma_{(0),\pm}^s+\sigma_{(0),\pm}^a+\sigma_{(0),\pm}^d\big).\label{equ:v-4}
\end{align}
Equating the $\varepsilon^k(k=0,1)$ terms in (\ref{BE-d}), we get
\ben
\na\cdot\vv^{(0)}_\pm=\na\cdot\vv^{(1)}_\pm=0.\label{equ:v-div}
\een
Equating the $\varepsilon^k(k=-2,-1,0)$ terms in (\ref{BE-Q}), we get
\begin{align}
&-f(\QQ^{(0)}_\pm)=\HH^{(-2)}_\pm=0,\label{equ:H-2-4}\\
&\pa_\tau\QQ^{(0)}_\pm=-f'(\QQ^{(0)}_\pm)\QQ^{(1)}_\pm=\HH^{(-1)}_\pm,\label{equ:H-1-4}\\
&\pa_\tau\QQ^{(1)}_\pm+\pa_t\QQ^{(0)}_\pm+\mathbf{v^{^{(0)}}_\pm}\cdot\nabla\QQ^{^{(0)}}_\pm+\QQ^{^{(0)}}_\pm\cdot\mathbf{\Omega^{^{(0)}}_\pm}
-\mathbf{\Omega^{^{(0)}}_\pm}\cdot\QQ^{^{(0)}}_\pm=\HH^{(0)}_\pm+\mathbf{S}_{\QQ^{(0)}_\pm}(\DD^{(0)}_\pm),\label{equ:Q-4}
\end{align}
where $\mathbf{\Omega^{^{(0)}}_\pm}=\f12\big(\nabla\mathbf{v^{(0)}_\pm}-(\nabla\mathbf{v^{(0)}_\pm})^{T}\big)$
and $\mathbf{D^{^{(0)}}_\pm}=\f12\big(\nabla\mathbf{v^{(0)}_\pm}+(\nabla\mathbf{v^{(0)}_\pm})^{T}\big)$.

By (\ref{equ:H-2-4}), (\ref{equ:H-1-4}), (\ref{equ:p-2-4}) and (\ref{equ:p-1-4}), we have
\beno
&&\sigma_{(-2),\pm}^s=0,\quad \sigma_{(-2),\pm}^a=0,\quad \na\mathbf{ p^{(-2)}_\pm}=0.
\eeno
Similar to analysis in Section 3.1, we know from (\ref{equ:H-2-4}) that
\begin{align}
&\QQ^{(0)}_\pm(\tau,t,x)=s_\pm\big(\nn(\tau,t,x)\nn(\tau,t,x)-\frac13\II\big)\label{equ:4.7}
\end{align}
for some $\nn(t,x)\in \BS$. The solvability of (\ref{equ:H-1-4}) to
find $\QQ^{(1)}_+$ requires that
\begin{align}
\pa_\tau\QQ^{(0)}_+\perp\text{ker}f'(\QQ^{(0)}_+).\label{equ:4.8}
\end{align}
By Proposition \ref{prop:kernel}, (\ref{equ:4.7}) and (\ref{equ:4.8}) require that for any $\nn^\perp\in \mathbf{V_n}$,
\begin{align}
&0=\big(\nn\nn_\tau+\nn_\tau\nn\big):\big(\nn\nn^{\bot}+\nn^{\bot}\nn\big)=2\nn_\tau\cdot\nn^\perp,\nonumber
\end{align}
which along with the fact that $\nn_\tau\cdot\nn=0$ and (\ref{equ:H-1-4}) implies that
\begin{align}
\nn_\tau=0,\qquad \pa_\tau\QQ_+^{(0)}=0,\qquad \HH^{(-1)}_+=0.\nonumber
\end{align}
Hence, we have
\beno
&&\sigma_{(-1),\pm}^s=0,\quad \sigma_{(-1),\pm}^a=0,
\eeno
and
\begin{align}
\pa_\tau\mathbf{v}^{(0)}_\pm=-\nabla\mathbf{p^{(-1)}}_\pm.
\end{align}
This implies that $\mathbf{p^{(-1)}}_\pm$ is harmonic by (\ref{equ:v-div}).

In order to derive the evolution equation of $(\vv^{(0)}, \QQ^{(0)})$  with respect to time scale $t$, we assume that as $\tau\rightarrow+\infty$,
\begin{align*}
&\mathbf{v}^{(0)}_\pm(\tau,t,x)\rightarrow\mathbf{v}^{(0)}_\pm(\infty,t,x),\quad \mathbf{p}^{(0)}_\pm(\tau,t,x)\rightarrow\mathbf{p}^{(0)}_\pm(\infty,t,x),
\\&\QQ^{(1)}_+(\tau,t,x)\rightarrow\QQ^{(1)}_+(\infty,t,x),\quad \pa_\tau\QQ^{(1)}_+(\tau,t,x)\rightarrow 0.
\end{align*}
Without confusing notations, we still use the same notations to denote the corresponding limits as $\tau\rightarrow+\infty$.

Recall that
\beno
\HH^{(0)}_+=\Delta\QQ^{(0)}_+-f'(\QQ^{(0)}_+)\QQ^{(2)}_+-\frac{1}{2}\big\langle f''(\QQ^{(0)}_+)\QQ^{(1)}_+,\QQ^{(1)}_\pm\big\rangle
\eeno
with $\QQ^{(1)}_+=\nn\nn^{\bot}+\nn^{\bot}\nn$ for $\nn^\bot\in \mathbf{V}_\nn$.
Then by (\ref{equ:Q-3.8}), the solvability of (\ref{equ:Q-4}) to find $\QQ^{(2)}_+$ requires that
\begin{align}
\pa_t\QQ^{(0)}_++\mathbf{v^{^{(0)}}_+}\cdot\nabla
\QQ^{^{(0)}}_++\QQ^{^{(0)}}_+\cdot\mathbf{\Omega^{^{(0)}}_+}
-\mathbf{\Omega^{^{(0)}}_+}\cdot\QQ^{^{(0)}}_+-\mathbf{S}_{\QQ^{(0)}_+}(\DD^{(0)}_+)-\Delta\QQ^{(0)}_+\perp \text{ker}f'(\QQ^{(0)}_+).\label{equ:Q-4.11}
\end{align}
Then by Lemma 3.2 in \cite{WZZ}, (\ref{equ:Q-4.11}) implies that $\nn(t,x)$ should satisfy
\beno
\nn\times\big(-s_+\Delta\nn+s_+\NN-\frac{\xi(2+s_+)}{3}\DD^{(0)}_+\cdot\nn\big)=0,
\eeno
where $\NN=\frac{\partial\nn}{\partial t}+\vv^{(0)}_+\cdot\nabla\nn-\BOm^{(0)}_+\cdot\nn$.
And by Lemma 3.3 in \cite{WZZ}, we know that
\begin{align*}
&\sigma_{(0),+}^s+\sigma_{(0),+}^a=\alpha_1(\nn\nn):\DD^{(0)}_+\nn\nn+\alpha_2\NN\nn+\alpha_3\nn\NN+\alpha_4\DD^{(0)}_++\alpha_5\DD^{(0)}_+\cdot\nn\nn\\
&\qquad\qquad\qquad\quad+\alpha_6\nn\nn\cdot\DD^{(0)}_+\triangleq \sigma^L,\\
&\sigma_{(0),+}^d=-2\na \nn\odot\na\nn\triangleq \sigma^E,
\end{align*}
where $\alpha_1,\cdots,\alpha_6$ are called the Leslie coefficients given by
\begin{align*}
 &\alpha_1=-\frac{2\xi^2s_+^2(3-2s_+)(1+2s_+)}{3},\quad \alpha_2=-s_+^2-\frac{\xi s_+(2+s_+)}{3}, \quad\alpha_3=s_+^2-\frac{\xi s_+(2+s_+)}{3},
\\&\alpha_4=1+\frac{4\xi^2(1-s_+)^2}{9},\quad \alpha_5=\frac{\xi^2s_+(4-s_+)}{3}+\frac{\xi s_+(2+s_+)}{3},\quad \alpha_6=\frac{\xi^2s_+(4-s_+)}{3}-\frac{\xi s_+(2+s_+)}{3}.
\end{align*}
Moreover, we have
\beno
\sigma_{(0),-}^s+\sigma_{(0),-}^a+\sigma_{(0),-}^d=\DD^{(0)}_-.
\eeno

In summary, in the region $\Om^+(t)$, $(\vv^{(0)}_+,\nn)$ satisfies
\begin{align*}
&\pa_t\mathbf{v}^{(0)}_++\mathbf{v^{(0)}_+\cdot\nabla v^{(0)}_+}=-\mathbf{\nabla
p^{(0)}_+}+\nabla\cdot\big(\sigma^L+\sigma^E\big),\\
&\na\cdot\vv^{(0)}_+=0,\\
&\nn\times\big(-\Delta\nn+\NN-\DD^{(0)}_+\cdot\nn\big)=0,
\end{align*}
which is the Ericksen-Leslie system introduced by Ericksen and Leslie \cite{Eri,Les}.
While in the region $\Om^-(t)$, $\vv^{(0)}_-$ satisfies the incompressible Navier-Stokes equations
\begin{align*}
&\pa_t\mathbf{v}^{(0)}_-+\mathbf{v^{(0)}_-\cdot\nabla v^{(0)}_-}=-\mathbf{\nabla
p^{(0)}_-}+\Delta \vv^{(0)}_-,\\
&\na\cdot\vv^{(0)}_-=0.
\end{align*}

\subsection{Inner expansion}

Using the fact that
\begin{align*}
\nabla\cdot\sigma_{\varepsilon}^d=-\nabla\big(\frac{1}{2}|\nabla\QQ^{\varepsilon}|^2+\varepsilon^{-2}
F_b(\QQ^{\varepsilon})\big)-\HH^{\varepsilon}:\nabla\QQ^{\varepsilon},
\end{align*}
the equation (\ref{BE-v}) can be rewritten as
\begin{align}\label{BE-v-m}
&\mathbf{v}^{\varepsilon}_t+\mathbf{v^{\varepsilon}\nabla
v^{\varepsilon}}-\Delta\mathbf{v^{\varepsilon}}+\mathbf{\nabla
q^{\varepsilon}}=\nabla\cdot(\sigma_{\varepsilon}^a-\mathbf{S}_{\QQ^{\varepsilon}}(\HH^{\varepsilon}))-\HH^{\varepsilon}:\nabla\QQ^{\varepsilon},
\end{align}
where $q^{\varepsilon}$ is the modified pressure defined by
\begin{align}
q^{\varepsilon}=p^{\varepsilon}+\frac{1}{2}|\nabla\QQ^{\varepsilon}|^2+\varepsilon^{-2}
F_b(\QQ^{\varepsilon}).\nonumber
\end{align}

In the transition region, we make the following expansion for $\big(\vv^\ve, \pp^\ve, \QQ^\ve\big)$:
\begin{align}
&\mathbf{v}^{\varepsilon}(t,x)=\widetilde{\mathbf{v}}^{(0)}(\tau,t,x,z)+\varepsilon\widetilde{\mathbf{v}}^{(1)}(\tau,t,x,z)+\cdots,\\&
\mathbf{p}^{\varepsilon}(t,x)=\varepsilon^{-2}\widetilde{\mathbf{p}}^{(-2)}(\tau,t,x,z)+\varepsilon^{-1}\widetilde{\mathbf{p}}^{(-1)}(\tau,t,x,z)+\cdots,\\&
\QQ^{\varepsilon}(t,x)=\widetilde{\QQ}^{(0)}(\tau,t,x,z)+\varepsilon\widetilde{\QQ}^{(1)}(\tau,t,x,z)+\cdots,
\end{align}
with $z=\f {\varphi(t,x)} \ve$. Then $\DD^\ve$ and $\mathbf{\Omega^{\varepsilon}}$ has the expansion
\begin{align*}
\DD^\ve=&\varepsilon^{-1}\frac{1}{2}\Big(\big(\nabla\varphi\mathbf{\widetilde{v}}^{(0)}_z\big)+\big(\nabla\varphi\mathbf{\widetilde{v}}_z^{(0)}\big)^{T}\Big)
+\frac{1}{2}\Big(\big(\nabla\varphi\mathbf{\widetilde{v}}^{(1)}_z\big)+\big(\nabla\varphi\mathbf{\widetilde{v}}_z^{(1)}\big)^{T}\\
&+\big(\nabla_x\mathbf{\widetilde{v}}^{(0)}\big)+\big(\nabla_x\mathbf{\widetilde{v}}^{(0)}\big)^{T}\Big)+\cdots
\triangleq \ve^{-1}\widetilde{\DD}^{(-1)}+\widetilde{\DD}^{(0)}+\cdots,\\
\mathbf{\Omega^{\varepsilon}}=&\varepsilon^{-1}\frac{1}{2}\Big(\big(\nabla\varphi\mathbf{\widetilde{v}}^{(0)}_z\big)-\big(\nabla\varphi\mathbf{\widetilde{v}}_z^{(0)}\big)^{T}\Big)
+\frac{1}{2}\Big(\big(\nabla\varphi\mathbf{\widetilde{v}}^{(1)}_z\big)-\big(\nabla\varphi\mathbf{\widetilde{v}}_z^{(1)}\big)^{T}\\
&+\big(\nabla_x\mathbf{\widetilde{v}}^{(0)}\big)-\big(\nabla_x\mathbf{\widetilde{v}}^{(0)}\big)^{T}\Big)+\cdots
\triangleq \ve^{-1}\widetilde{\mathbf\Omega}^{(-1)}+\widetilde{\mathbf\Omega}^{(0)}+\cdots,
\end{align*}
The molecular field $\HH^\ve$ has the expansion
\begin{align*}
\HH^\ve=&\ve^{-2}\Big(\widetilde{\QQ}^{(0)}_{zz}-f
(\widetilde{\QQ}^{(0)})\Big)+\ve^{-1}\Big(\widetilde{\QQ}^{(1)}_{zz}-f'(\widetilde{\QQ}^{(0)})\widetilde{\QQ}^{(1)}\\
&\quad+\widetilde{\QQ}^{(0)}_{z}\Delta\varphi+2\nabla\varphi\cdot\nabla_x\widetilde{\QQ}^{(0)}_{z}\Big)+\Big(\widetilde{\QQ}^{(2)}_{zz}-f'(\widetilde{\QQ}^{(0)})\widetilde{\QQ}^{(2)}
-\frac{1}{2}\big\langle f''(\widetilde{\QQ}^{(0)})\widetilde{\QQ}^{(1)},\widetilde{\QQ}^{(1)}\big\rangle\\
&\quad+\widetilde{\QQ}^{(1)}_z\Delta\varphi+2\na\varphi\cdot\na_x\widetilde{\QQ}^{(1)}_z+\Delta_x\widetilde{\QQ}^{(0)}\Big)+\cdots\\
\triangleq& \ve^{-2}\widetilde{\HH}^{(-2)}+\ve^{-1}\widetilde{\HH}^{(-1)}+\widetilde{\HH}^{(0)}+\cdots.
\end{align*}
The anti-symmetry viscous stress $\sigma_\ve^a$  has the expansion
\begin{align*}
\sigma_\ve^a=&\ve^{-2}\Big(\widetilde{\QQ}^{(0)}\cdot\widetilde{\HH}^{(-2)}-\widetilde{\HH}^{(-2)}\cdot\widetilde{\QQ}^{(0)}\Big)
+\ve^{-1}\Big(\widetilde{\QQ}^{(0)}\cdot\widetilde{\HH}^{(-1)}-\widetilde{\HH}^{(-1)}\cdot\widetilde{\QQ}^{(0)}\\
&\qquad+\widetilde{\QQ}^{(1)}\cdot\widetilde{\HH}^{(-2)}-\widetilde{\HH}^{(-2)}\cdot\widetilde{\QQ}^{(1)}\Big)+\Big(\widetilde{\QQ}^{(0)}\cdot\widetilde{\HH}^{(0)}-\widetilde{\HH}^{(0)}\cdot\widetilde{\QQ}^{(0)}\Big)+\cdots\\
\triangleq& \ve^{-2}\widetilde{\sigma}^a_{(-2)}+\ve^{-1}\widetilde{\sigma}^a_{(-1)}+\widetilde{\sigma}^a_{(0)}\cdots.
\end{align*}
The modified pressure $\widetilde{q}^\ve$ has the expansion
\begin{align}
\widetilde{\mathbf{q}}^{\varepsilon}=&\varepsilon^{-2}\Big(\frac{1}{2}|\widetilde{\QQ}_{z}^{^{(0)}}|^2+F_b(\widetilde{\QQ}^{(0)})
+\widetilde{\mathbf{p}}^{(-2)}\Big)+\varepsilon^{-1}\Big(\widetilde{\QQ}_{z}^{^{(0)}}:\widetilde{\QQ}_{z}^{^{(1)}}
\nonumber\\&+f(\widetilde{\QQ}^{(0)})\widetilde{\QQ}^{(1)}+\nabla\varphi\cdot\nabla\widetilde{\QQ}^{(0)}:\widetilde{\QQ}_z^{(0)}+\widetilde{\mathbf{p}}^{(-1)}\Big)
\nonumber\\&+\Big(\widetilde{\QQ}_{z}^{^{(0)}}:\widetilde{\QQ}_{z}^{^{(2)}}+\frac{1}{2}|\widetilde{\QQ}_{z}^{^{(1)}}|^2
+f(\widetilde{\QQ}^{(0)})\widetilde{\QQ}^{(2)}+\frac{1}{2}\langle f'(\widetilde{\QQ}^{(0)})\widetilde{\QQ}^{(1)},\widetilde{\QQ}^{(1)}\rangle
\nonumber\\&\quad+\nabla\varphi\cdot\nabla\widetilde{\QQ}^{(0)}:\widetilde{\QQ}_z^{(1)}+\nabla\varphi\cdot\nabla\widetilde{\QQ}^{(1)}:\widetilde{\QQ}_z^{(0)}
+\frac{1}{2}\nabla_{x}\widetilde{\QQ}^{(0)}:\nabla_{x}\widetilde{\QQ}^{(0)}+\widetilde{\mathbf{p}}^{(0)}\Big)+\cdots
\nonumber\\\triangleq &\varepsilon^{-2}\widetilde{\mathbf{q}}^{(-2)}+\varepsilon^{-1}\widetilde{\mathbf{q}}^{(-1)}+\widetilde{\mathbf{q}}^{(0)}+\cdots.\nonumber
\end{align}

Now we plugg these expansions into (\ref{BE-v-m}), then equate $\ve^k(k=-3,-2,-1)$ terms to obtain
\begin{align}
&\partial_z\big(\widetilde{\sigma}_{-2}^a-\mathbf{S}_{\widetilde{\QQ}^{(0)}}(\widetilde{\mathbf{H}}^{(-2)})\big)\cdot\nabla\varphi
=\widetilde{\HH}^{(-2)}:\nabla\varphi\widetilde{\QQ}^{(0)}_z+\nabla\varphi\widetilde{\mathbf{q}}_z^{(-2)},\label{equ:v-exp-2}\\
&\frac{\nu}{2}\widetilde{\mathbf{v}}^{(0)}_{zz}
=\varphi_\tau\widetilde{\mathbf{v}}^{(0)}_{z}-\partial_z\big(\widetilde{\sigma}_{(-1)}^a-\mathbf{S}_{\widetilde{\QQ}^{(0)}}(\widetilde{\mathbf{H}}^{(-1)})
-\mathbf{S}_{\widetilde{\QQ}^{(1)}}(\widetilde{\mathbf{H}}^{(-2)})\big)\cdot\nabla\varphi\nonumber\\
&\qquad\qquad\qquad-\nabla\cdot\big(\widetilde{\sigma}_{(-2)}^a-\mathbf{S}_{\widetilde{\QQ}^{(0)}}(\widetilde{\mathbf{H}}^{(-2)})\big)
+\nabla\varphi\widetilde{\mathbf{q}}_z^{(-1)}+\nabla_x\widetilde{\mathbf{q}}^{(-2)}
\nonumber\\&\qquad\qquad\qquad+\big(\widetilde{\HH}^{(-2)}:\widetilde{\QQ}^{(1)}_z+\widetilde{\HH}^{(-1)}:\widetilde{\QQ}^{(0)}_z\big)\nabla\varphi
+\widetilde{\HH}^{(-2)}:\nabla_x\widetilde{\QQ}^{(0)},\label{equ:v-exp-1}
\end{align}
and
\begin{align}
\widetilde{\mathbf{v}}^{(1)}_{zz}
=&\varphi_\tau\widetilde{\mathbf{v}}^{(1)}_{z}-\partial_z\big(\widetilde{\sigma}_{(0)}^a-\mathbf{S}_{\widetilde{\QQ}^{(2)}}(\widetilde{\mathbf{H}}^{(-2)})
-\mathbf{S}_{\widetilde{\QQ}^{(1)}}(\widetilde{\mathbf{H}}^{(-1)})-\mathbf{S}_{\widetilde{\QQ}^{(0)}}(\widetilde{\mathbf{H}}^{(0)})\big)\cdot\nabla\varphi
\nonumber\\&-\nabla_x\cdot\big(\widetilde{\sigma}_{(-1)}^a-\mathbf{S}_{\widetilde{\QQ}^{(1)}}(\widetilde{\mathbf{H}}^{(-2)})
-\mathbf{S}_{\widetilde{\QQ}^{(0)}}(\widetilde{\mathbf{H}}^{(-1)})\big)+\nabla\varphi\widetilde{\mathbf{q}}_z^{(0)}+\nabla_x\widetilde{\mathbf{q}}^{(-1)}
\nonumber\\&+\big(\varphi_t\widetilde{\mathbf{v}}^{(0)}_z+\widetilde{\mathbf{v}}^{(0)}_\tau+\widetilde{\mathbf{v}}^{(0)}\cdot\nabla\varphi\widetilde{\mathbf{v}}^{(0)}_z
-\widetilde{\mathbf{v}}_{z}^{(0)}\Delta\varphi-2\nabla\varphi\cdot\nabla_x\widetilde{\mathbf{v}}_{z}^{(0)}\big)
\nonumber\\&+\big(\widetilde{\HH}^{(-2)}:\widetilde{\QQ}^{(2)}_z+\widetilde{\HH}^{(-1)}:\widetilde{\QQ}^{(1)}_z
+\widetilde{\HH}^{(0)}:\widetilde{\QQ}^{(0)}_z\big)\nabla\varphi
\nonumber\\&+\big(\widetilde{\HH}^{(-2)}:\nabla_x\widetilde{\QQ}^{(1)}+\widetilde{\HH}^{(-1)}:\nabla_x\widetilde{\QQ}^{(0)}\big).\label{equ:v-exp-0}
\end{align}
Equating $\varepsilon^k(k=-1,0)$ terms  in (\ref{BE-d}), we obtain
\begin{align}
&\nabla\varphi\cdot\mathbf{\widetilde{v}}_z^{(0)}=0,\label{equ:d-exp-1}\\&
\nabla\varphi\cdot\mathbf{\widetilde{v}}_z^{(1)}+\nabla_x\cdot
\widetilde{\mathbf{v}}^{(0)}=0.\label{equ:d-exp-0}
\end{align}
Equating $\varepsilon^k(k=-2,-1)$ terms in (\ref{BE-Q}),  we obtain
\begin{align}
&-\varphi_\tau\widetilde{\QQ}^{(0)}_z+\widetilde{\QQ}_{zz}^{^{(0)}}-f(\widetilde{\QQ}^{(0)})=0,\label{equ:Q-exp-2}
\\&\widetilde{\HH}^{(-1)}+\varphi_\tau\widetilde{\QQ}^{(1)}_z=\big(\nabla\varphi\cdot\mathbf{\widetilde{v}^{(0)}}+\varphi_t\big)\widetilde{\QQ}^{(0)}_{z}\nonumber\\
&\qquad\qquad-\widetilde{\QQ}^{(0)}_\tau-\mathbf{\widetilde{\Omega}}^{(-1)}\cdot\widetilde{\QQ}^{(0)}+\widetilde{\QQ}^{(0)}\cdot\mathbf{\widetilde{\Omega}}^{(-1)}
-\mathbf{S}_{\widetilde{\QQ}^{(0)}}\cdot(\widetilde{\mathbf{D}}^{(-1)}).\label{equ:Q-exp-1}
\end{align}

Multiplying (\ref{equ:Q-exp-2}) by $\widetilde{\QQ}^{(0)}$ and integrating from
$z=-\infty$ to $z=+\infty$,  we obtain
\begin{align}
\varphi_\tau\int_{-\infty}^{\infty}|\widetilde{\QQ}_z^{(0)}|^2dz&=F_b(\widetilde{\QQ}^{(0)}(\tau,t,x,-\infty))-F_b(\widetilde{\QQ}^{(0)}(\tau,t,x,+\infty))=0.
\end{align}
Hence,
\begin{align}
\varphi_\tau=0, \quad -\widetilde{\QQ}_{zz}^{(0)}+f(\widetilde{\QQ}^{(0)})=0,\label{equ:4.26}
\end{align}
and  $\widetilde{\HH}^{(-2)}=0$ by (\ref{equ:Q-exp-2}). Hence,
\ben
\widetilde{\sigma}^a_{-2}=0,\quad \widetilde{\sigma}^a_{-1}=\widetilde{\QQ}^{(0)}\cdot\widetilde{\HH}^{(-1)}-\widetilde{\HH}^{(-1)}\cdot\widetilde{\QQ}^{(0)},\label{equ:4.24}
\een
and by (\ref{equ:v-exp-2}),
\ben
 \widetilde{\mathbf{q}}_z^{(-2)}=\pa_z\Big(\frac{1}{2}|\widetilde{\QQ}_{z}^{^{(0)}}|^2+F_b(\widetilde{\QQ}^{(0)})
+\widetilde{\mathbf{p}}^{(-2)}\Big)=0.\label{equ:q-2}
\een
Then the equation (\ref{equ:v-exp-1}) is reduced to
\begin{align}
&\widetilde{\mathbf{v}}^{(0)}_{zz}
=-\partial_z\big(\widetilde{\sigma}_{(-1)}^a-\mathbf{S}_{\widetilde{\QQ}^{(0)}}(\widetilde{\mathbf{H}}^{(-1)})\big)
\cdot\nabla\varphi+\nabla\varphi\widetilde{\mathbf{q}}_z^{(-1)}+\nabla_x\widetilde{\mathbf{q}}^{(-2)}
\nonumber\\&\qquad\qquad\qquad+\widetilde{\HH}^{(-1)}:\widetilde{\QQ}^{(0)}_z\nabla\varphi.\label{equ:v-exp-new1}
\end{align}
While, the equation (\ref{equ:v-exp-0}) is reduced to
\begin{align}
\widetilde{\mathbf{v}}^{(1)}_{zz}
=&-\partial_z\big(\widetilde{\sigma}_{(0)}^a-\mathbf{S}_{\widetilde{\QQ}^{(1)}}(\widetilde{\mathbf{H}}^{(-1)})-\mathbf{S}_{\widetilde{\QQ}^{(0)}}(\widetilde{\mathbf{H}}^{(0)})\big)\cdot\nabla\varphi
\nonumber\\&-\nabla_x\cdot\big(\widetilde{\sigma}_{(-1)}^a-\mathbf{S}_{\widetilde{\QQ}^{(0)}}(\widetilde{\mathbf{H}}^{(-1)})\big)
+\nabla\varphi\widetilde{\mathbf{q}}_z^{(0)}+\nabla_x\widetilde{\mathbf{q}}^{(-1)}
\nonumber\\&+\big(\varphi_t\widetilde{\mathbf{v}}^{(0)}_z+\widetilde{\mathbf{v}}^{(0)}_\tau+\widetilde{\mathbf{v}}^{(0)}\cdot\nabla\varphi\widetilde{\mathbf{v}}^{(0)}_z
-\widetilde{\mathbf{v}}_{z}^{(0)}\Delta\varphi-2\nabla\varphi\cdot\nabla_x\widetilde{\mathbf{v}}_{z}^{(0)}\big)
\nonumber\\&+\big(\widetilde{\HH}^{(-1)}:\widetilde{\QQ}^{(1)}_z
+\widetilde{\HH}^{(0)}:\widetilde{\QQ}^{(0)}_z\big)\nabla\varphi+\widetilde{\HH}^{(-1)}:\nabla_x\widetilde{\QQ}^{(0)}\big).\label{equ:v-exp-new0}
\end{align}

\subsection{Evolution of the sharp interface}

By (\ref{equ:4.26}), it is reasonable to assume that $\widetilde{\QQ}^{(0)}(\tau,t,x,z)=s(z)(\nn(t, x)\nn(t,x)-\frac{1}{3}\II)$, hence $\pa_\tau\widetilde{\QQ}^{(0)}=0$.
Thanks to the definition of $\widetilde{\HH}^{(-1)}$, the equation (\ref{equ:Q-exp-1}) can be rewritten as
\begin{align}
&\widetilde{\QQ}^{(1)}_{zz}-f'(\widetilde{\QQ}^{(0)})\widetilde{\QQ}^{(1)}\nonumber\\
&=\big(\nabla\varphi\cdot\mathbf{\widetilde{v}^{(0)}}+\varphi_t-\Delta\varphi-2\na\varphi\cdot\na_x\big)\widetilde{\QQ}^{(0)}_{z}
-\mathbf{\widetilde{\Omega}}^{(-1)}\cdot\widetilde{\QQ}^{(0)}+\widetilde{\QQ}^{(0)}\cdot\mathbf{\widetilde{\Omega}}^{(-1)}
-\mathbf{S}_{\widetilde{\QQ}^{(0)}}\cdot(\widetilde{\mathbf{D}}^{(-1)}).\label{equ:Q-exp-new1}
\end{align}
Multiplying (\ref{equ:Q-exp-new1}) by $\widetilde{\QQ}_{z}^{(0)}$ and integrating from
$z=-\infty$ to $z=+\infty$, we obtain
\begin{align}
\varphi_t&=\Delta\varphi-\nabla\varphi\cdot\mathbf{v}^{(0)}+\frac{\nabla\varphi\cdot\nabla\int_{-\infty}^{+\infty}|\widetilde{\QQ}_{z}^{(0)}|^2dz}{\int_{-\infty}^{+\infty}|\widetilde{\QQ}_{z}^{(0)}|^2dz}
\nonumber\\&\qquad+\bigg(\int_{-\infty}^{+\infty}|\widetilde{\QQ}_{z}^{(0)}|^2dz\bigg)^{-1}
\int_{-\infty}^{\infty}\Big(\mathbf{\widetilde{\Omega}}^{(-1)}\cdot\widetilde{\QQ}^{(0)}-\widetilde{\QQ}^{(0)}\cdot\mathbf{\widetilde{\Omega}}^{(-1)}
+\mathbf{S}_{\widetilde{\QQ}^{(0)}}(\widetilde{\mathbf{D}}^{(-1)})\Big):\widetilde{\QQ}_{z}^{(0)} dz
\nonumber\\&=\Delta\varphi-\nabla\varphi\cdot\mathbf{v}^{(0)}+\frac{\nabla\varphi\cdot\nabla\int_{-\infty}^{+\infty}|\widetilde{\QQ}_{z}^{(0)}|^2dz}{\int_{-\infty}^{+\infty}|\widetilde{\QQ}_{z}^{(0)}|^2dz}
\nonumber\\&\qquad+\bigg(\int_{-\infty}^{+\infty}|\widetilde{\QQ}_{z}^{(0)}|^2dz\bigg)^{-1}
\int_{-\infty}^{\infty}\mathbf{\widetilde{\Omega}}^{(-1)}:\big(\widetilde{\QQ}_{z}^{(0)}\widetilde{\QQ}^{(0)}-\widetilde{\QQ}^{(0)}\widetilde{\QQ}_{z}^{(0)}\big) dz
\nonumber\\&\qquad+\xi\bigg(\int_{-\infty}^{+\infty}|\widetilde{\QQ}_{z}^{(0)}|^2dz\bigg)^{-1}
\int_{-\infty}^{\infty}\widetilde{\mathbf{D}}^{(-1)}:\big(\widetilde{\QQ}_{z}^{(0)}\widetilde{\QQ}^{(0)}+\widetilde{\QQ}^{(0)}\widetilde{\QQ}_{z}^{(0)}\big) dz
\nonumber\\&\qquad+\xi\bigg(\int_{-\infty}^{+\infty}|\widetilde{\QQ}_{z}^{(0)}|^2dz\bigg)^{-1}
\int_{-\infty}^{\infty}\Big(\widetilde{\mathbf{D}}^{(-1)}:\frac{2}{3}\widetilde{\QQ}_{z}^{(0)}
-2\big(\widetilde{\QQ}^{(0)}:\widetilde{\QQ}_{z}^{(0)}\big)\big(\widetilde{\mathbf{D}}^{(-1)}:\widetilde{\QQ}^{(0)}\big)\Big)dz.\label{equ:MCF}
\end{align}
Here we used the fact that $\nabla\varphi\cdot\mathbf{\widetilde{v}^{(0)}}=\nabla\varphi\cdot\mathbf{{v}^{(0)}}$ by (\ref{equ:d-exp-1}).

As $\int_{-\infty}^{+\infty}|\widetilde{\QQ}_{z}^{(0)}|^2dz$ is independent of $(t,x)$, the equation (\ref{equ:MCF}) is reduced to
\begin{align}
\varphi_t&=\Delta\varphi-\nabla\varphi\cdot\mathbf{v}^{(0)}
\nonumber\\&\quad+\frac{\xi}{2}\bigg(\int_{-\infty}^{+\infty}|s'(z)|^2dz\bigg)^{-1}\int_{-\infty}^{\infty}s'(1+s-2s^2)\widetilde{\mathbf{D}}^{(-1)}:(\nn\nn)dz.\label{equ:MCF-f}
\end{align}

\subsection{Jump conditions on the sharp interface}
In this subsection, we will derive the jump condition of the velocity and the pressure.
For this end, we assume that as $z\rightarrow \pm \infty$,
\begin{align*}
&\widetilde{\mathbf{v}}^{(k)}(\tau,t,x,z)\rightarrow
\mathbf{v}^{(k)}_\pm(\tau,t,x)\quad k=0,1,\\
&\widetilde{\QQ}^{(k)}(\tau,x,z)\rightarrow \QQ^{(k)}_\pm(\tau,t,x)\quad k=0,1,\\
&\widetilde{\mathbf{p}}^{(k)}(\tau,t,x,z)\rightarrow \mathbf{p}^{(k)}_\pm(\tau,t,x)\quad k=-2,-1,0.
\end{align*}
We also assume that $\widetilde{\QQ}^{(0)}(\tau,t,x,z)=s(z)(\nn(t, x)\nn(t,x)-\frac{1}{3}\II)$.
In the sequel, $\nu$ is the unit normal to the sharp interface $\Gamma(t)$ and  $[\cdot]$ denotes the jump
across $\Gamma(t)$, i.e.,
\beno
[f]=f_+-f_-\big|_{\Gamma(t)}.
\eeno

First of all, integrating (\ref{equ:d-exp-1}) from $z=-\infty$ to $z=+\infty$, we derive the jump condition of the velocity on $\Gamma(t)$:
\begin{align}
\big[\mathbf{v}^{(0)}\big]\cdot \nu=0\quad\text{on}\quad \Gamma(t).
\end{align}

Integrating (\ref{equ:q-2}) from $z=-\infty$ to $z=+\infty$, we obtain
\begin{equation}
\big[\mathbf{p}^{(-2)}\big]=-\big[F_b(\QQ^{(0)})\big]=0\quad\text{on}\quad \Gamma(t).
\end{equation}

Note that $\na_x\widetilde{\mathbf{q}}^{(-2)}=\na_x\mathbf{p}^{(-2)}=0$. Multiplying (\ref{equ:v-exp-new1}) by $\nabla\varphi$,  we get by (\ref{equ:d-exp-1}) that
\begin{align}\label{equ:q-4.30}
&\widetilde{\mathbf{q}}_z^{(-1)}=\partial_z\Big\langle \widetilde{\sigma}_{(-1)}^a-\mathbf{S}_{\widetilde{\QQ}^{(0)}}(\widetilde{\mathbf{H}}^{(-1)}),
{\nabla\varphi}\otimes{\nabla\varphi}\Big\rangle-\widetilde{\HH}^{(-1)}:\widetilde{\QQ}^{(0)}_z.
\end{align}
Then we get by integrating (\ref{equ:q-4.30}) from $z=-\infty$ to $z=+\infty$ that
\begin{align}
\big[\mathbf{q}^{(-1)}\big]=&-\int_{-\infty}^\infty\widetilde{\HH}^{(-1)}:\widetilde{\QQ}^{(0)}_zdz\nonumber\\
=&-\big(\Delta\varphi+\nabla\varphi\cdot\nabla\big)\int_{-\infty}^{+\infty}\big|\widetilde{\QQ}_z^{(0)}\big|^2dz\nonumber\\
=&-\frac{2}{3}\Delta\varphi\int_{-\infty}^{+\infty}\big|s'(z)\big|^2dz,\label{equ:jump-q-1}
\end{align}
here we used the fact that
\begin{align}
&\int_{-\infty}^{+\infty}\big\langle\widetilde{\QQ}_{zz}^{(1)}-f'(\widetilde{\QQ}^{(0)})\widetilde{\QQ}^{(1)},\widetilde{\QQ}_{z}^{(0)}\big\rangle dz
\nonumber\\&=\int_{-\infty}^{+\infty}\big\langle\widetilde{\QQ}_{zzz}^{(0)}-f'(\widetilde{\QQ}^{(0)})\widetilde{\QQ}^{(0)}_z,\widetilde{\QQ}^{(1)}\big\rangle dz=0
\quad \textrm{by }(\ref{equ:4.26}).
\nonumber
\end{align}
By (\ref{equ:4.26}) again, we have
\begin{align}
\widetilde{\mathbf{q}}^{(-1)}&=\widetilde{\QQ}_{z}^{^{(0)}}:\widetilde{\QQ}_{z}^{^{(1)}}+f(\widetilde{\QQ}^{(0)})\widetilde{\QQ}^{(1)}+\nabla\varphi\cdot\nabla\widetilde{\QQ}^{(0)}:\widetilde{\QQ}_z^{(0)}+\widetilde{\mathbf{p}}^{(-1)}
\nonumber\\&=\partial_z\big(\widetilde{\QQ}_{z}^{^{(0)}}:\widetilde{\QQ}^{(1)}\big)
+\nabla\varphi\cdot\nabla\widetilde{\QQ}^{(0)}:\widetilde{\QQ}_z^{(0)}+\widetilde{\mathbf{p}}^{(-1)},\label{equ:q-1/p}
\end{align}
which along with (\ref{equ:jump-q-1}) implies that
\begin{equation}
\big[\mathbf{p}^{(-1)}\big]=\big[\mathbf{q}^{(-1)}\big]
=-\frac{2}{3}\Delta\varphi\int_{-\infty}^{+\infty}\big|s'(z)\big|^2dz.\label{equ:jump-p-1}
\end{equation}

Next we derive the jump condition of $\mathbf{p}^{(0)}_\pm$. We get by (\ref{equ:d-exp-1})  that
\beno
&&\big(\varphi_t\widetilde{\mathbf{v}}^{(0)}_z+\widetilde{\mathbf{v}}^{(0)}_\tau+\widetilde{\mathbf{v}}^{(0)}\cdot\nabla\varphi\widetilde{\mathbf{v}}^{(0)}_z
-\widetilde{\mathbf{v}}_{z}^{(0)}\Delta\varphi-2\nabla\varphi\cdot\nabla_x\widetilde{\mathbf{v}}_{z}^{(0)}\big)\cdot\na\varphi\\
&&=\widetilde{\mathbf{v}}^{(0)}_\tau\cdot\na\varphi-\big\langle2\widetilde{\mathbf{D}}^{(0)},{\nabla\varphi}\otimes{\nabla\varphi}\big\rangle,
\eeno
and by (\ref{equ:d-exp-0}), we have
\beno
\widetilde{\mathbf{v}}^{(1)}_{zz}\cdot\na \varphi=-\big(\nabla_x\cdot\widetilde{\mathbf{v}}^{(0)}\big)_{z}.
\eeno
Then multiply (\ref{equ:v-exp-new0}) by $\nabla\varphi$ to obtain
\begin{align}
\widetilde{\mathbf{q}}_z^{(0)}=&-\big(\nabla_x\widetilde{\mathbf{q}}^{(-1)}-\nabla\mathbf{p^{(-1)}}_\pm\big)\cdot{\nabla\varphi}+\big(\nabla_x\cdot\widetilde{\mathbf{v}}^{(0)}\big)_{z}
\nonumber\\&+\partial_z\Big\langle\widetilde{\sigma}_{(0)}^a+2\widetilde{\mathbf{D}}^{(0)}-\mathbf{S}_{\widetilde{\QQ}^{(1)}}(\widetilde{\mathbf{H}}^{(-1)})-\mathbf{S}_{\widetilde{\QQ}^{(0)}}(\widetilde{\mathbf{H}}_{0}),{\nabla\varphi}\otimes{\nabla\varphi}\Big\rangle
\nonumber\\&+\Big(\nabla_x\cdot\big(\widetilde{\sigma}_{(-1)}^a-\mathbf{S}_{\widetilde{\QQ}^{(0)}}(\widetilde{\mathbf{H}}^{(-1)})\big)\Big)\cdot{\nabla\varphi}
\nonumber\\&-\big(\widetilde{\HH}^{(-1)}:\widetilde{\QQ}^{(1)}_z+\widetilde{\HH}^{(0)}:\widetilde{\QQ}^{(0)}_z\big)
-\big(\widetilde{\HH}^{(-1)}:\nabla_x\widetilde{\QQ}^{(0)}\big)\cdot{\nabla\varphi}.\label{equ:qz-1}
\end{align}
Here we used the fact that
\beno
\na\varphi\cdot\widetilde{\mathbf{v}}^{(0)}_\tau=\big(\na\varphi\cdot{\mathbf{v}}^{(0)}_\pm\big)_\tau=-\nabla\mathbf{p^{(-1)}}_\pm\cdot \na\varphi.
\eeno
Integrating  (\ref{equ:qz-1}) from $z=-\infty$ to $z=+\infty$, we get by (\ref{equ:4.24}) that
\begin{align}
&\big[\mathbf{q}^{(0)}\big]-\bigg[\Big\langle\widetilde{\sigma}_{(0)}^a+2\widetilde{\mathbf{D}}^{(0)}
-\mathbf{S}_{\widetilde{\QQ}^{(1)}}(\widetilde{\mathbf{H}}^{(-1)})-\mathbf{S}_{\widetilde{\QQ}^{(0)}}(\widetilde{\mathbf{H}}_{0}),{\nabla\varphi}\otimes{\nabla\varphi}\Big\rangle\bigg]
\nonumber\\&=-\int_{-\infty}^{+\infty}\big(\nabla_x\widetilde{\mathbf{q}}^{(-1)}-\nabla\mathbf{p^{(-1)}}_\pm\big)\cdot{\nabla\varphi}dz
\nonumber\\&\quad+\int_{-\infty}^{+\infty}\Big(\nabla_x\cdot\big(\widetilde{\QQ}^{(0)}\widetilde{\HH}^{(-1)}-\widetilde{\HH}^{(-1)}\widetilde{\QQ}^{(0)}
-\mathbf{S}_{\widetilde{\QQ}^{(0)}}(\widetilde{\HH}^{(-1)})\big)\Big)\cdot{\nabla\varphi}dz
\nonumber\\&\quad-\int_{-\infty}^{+\infty}\big(\widetilde{\HH}^{(-1)}:\widetilde{\QQ}^{(1)}_z
+\widetilde{\HH}^{(0)}:\widetilde{\QQ}^{(0)}_z\big)+\big(\widetilde{\HH}^{(-1)}:\nabla_x\widetilde{\QQ}^{(0)}\big)\cdot{\nabla\varphi}dz.\label{equ:q-41}
\end{align}
Thanks to the definition of $\widetilde{\mathbf{q}}^{(0)}$, we get by (\ref{equ:H-2-4}) and (\ref{equ:H-1-4}) that
\begin{align}
&\big[\mathbf{q}^{(0)}\big]=\int_{-\infty}^{+\infty}\nabla_x\widetilde{\QQ}^{(0)}:\nabla_x\widetilde{\QQ}^{(0)}_zdz+\big[\mathbf{p}^{(0)}\big],\label{equ:q-42}\\
&\Big[\widetilde{\sigma}_{0}^a-\mathbf{S}_{\widetilde{\QQ}^{(1)}}(\widetilde{\mathbf{H}}^{(-1)})\Big]
=\Big[\QQ^{(0)}\cdot\HH^{(0)}-\HH^{(0)}\cdot\QQ^{(0)}\Big].
\end{align}
And by (\ref{equ:d-exp-1}), we find
\begin{align}
\Big[\big\langle{\mathbf{{D}}}^{(0)},{\nabla\varphi}\otimes{\nabla\varphi}\big\rangle\Big]
&=\frac{1}{2}\int_{-\infty}^{+\infty}\big\langle\nabla\widetilde{\mathbf{v}}_z^{(0)}+(\nabla\widetilde{\mathbf{v}}_z^{(0)})^{T},{\nabla\varphi}\otimes{\nabla\varphi}\big\rangle dz
\nonumber\\&=\int_{-\infty}^{+\infty}\big\langle\nabla\widetilde{\mathbf{v}}_z^{(0)},{\nabla\varphi}\otimes{\nabla\varphi}\big\rangle dz
\nonumber\\&=\int_{-\infty}^{+\infty}\Big(\nabla\big(\nabla\varphi\cdot\widetilde{\mathbf{v}}_z^{(0)}\big)\nabla\varphi-\frac{1}{2}\widetilde{\mathbf{v}}_z^{(0)}\cdot\na |\nabla\varphi|^2\Big)dz
=0.
\end{align}
Thanks to the definition of  $\sigma^d$, we have
\begin{align}
&\int_{-\infty}^{+\infty}\Big(\big(\widetilde{\QQ}^{(0)}_z\Delta\varphi+2\na_x\widetilde{\QQ}^{(0)}_z\na\varphi\big):\nabla_x\widetilde{\QQ}^{(0)}\Big)
\cdot{\nabla\varphi}dz\nonumber\\
&=\int_{-\infty}^{+\infty}\Big(\widetilde{\QQ}^{(0)}_z\Delta\varphi:\nabla_x\widetilde{\QQ}^{(0)}\Big)\cdot{\nabla\varphi}dz
-\Big\langle\sigma^d(\QQ^{(0)},\QQ^{(0)}),{\nabla\varphi}\otimes{\nabla\varphi}\Big\rangle\nonumber\\
&=-\Big\langle\sigma^d(\QQ^{(0)},\QQ^{(0)}),{\nabla\varphi}\otimes{\nabla\varphi}\Big\rangle.\label{equ:4.41}
\end{align}
Let
\beno
&&\widetilde{\HH}^{(-1)}_1=\widetilde{\QQ}^{(1)}_{zz}-f'(\widetilde{\QQ}^{(0)})\widetilde{\QQ}^{(1)},\\
&&\widetilde{\HH}^{(0)}_1=\widetilde{\QQ}^{(2)}_{zz}-f'(\widetilde{\QQ}^{(0)})\widetilde{\QQ}^{(2)}
-\frac{1}{2}\big\langle f''(\widetilde{\QQ}^{(0)})\widetilde{\QQ}^{(1)},\widetilde{\QQ}^{(1)}\big\rangle.
\eeno
We find that
\begin{align}
&\int_{-\infty}^{+\infty}\big(\widetilde{\HH}_{1}^{(0)}:\widetilde{\QQ}^{(0)}_z+\widetilde{\HH}^{(-1)}_{1}:\widetilde{\QQ}^{(1)}_z\big)dz=0\label{equ:q-44}
\end{align}
by the following two facts
\begin{align}
&\int_{-\infty}^{+\infty}\big(f'(\widetilde{\QQ}^{(0)})\widetilde{\QQ}^{(1)}\big):\widetilde{\QQ}^{(1)}_z dz+\frac{1}{2}\int_{-\infty}^{+\infty}\big\langle f''(\widetilde{\QQ}^{(0)})\widetilde{\QQ}^{(1)},\widetilde{\QQ}^{(1)}\big\rangle:\widetilde{\QQ}^{(0)}_zdz
\nonumber\\&=\int_{-\infty}^{+\infty}\partial_z\Big(\frac{a}{2}|\widetilde{\QQ}^{(1)}|^2-\frac{b}{2}\big(\widetilde{\QQ}^{(1)}\cdot\widetilde{\QQ}^{(0)}
+\widetilde{\QQ}^{(0)}\cdot\widetilde{\QQ}^{(1)}\big):\widetilde{\QQ}^{(1)}+\frac{c}{2}|\widetilde{\QQ}^{(1)}|^2|\widetilde{\QQ}^{(0)}|^2
+c\big(\widetilde{\QQ}^{(1)}:\widetilde{\QQ}^{(0)}\big)^2\Big)dz
\nonumber\\&=\frac{1}{2}\int_{-\infty}^{+\infty}\partial_z\big(\big(f'(\widetilde{\QQ}^{(0)})\widetilde{\QQ}^{(1)}\big):\widetilde{\QQ}^{(1)}\big)dz
=\frac{1}{2}\Big[\big(f'({\QQ}^{(0)}){\QQ}^{(1)}\big):{\QQ}^{(1)}\Big]=0,\nonumber
\end{align}
and
\begin{align}
&\int_{-\infty}^{+\infty}\big(\widetilde{\QQ}_{zz}^{(2)}-f'(\widetilde{\QQ}^{(0)})\widetilde{\QQ}^{(2)}\big):\widetilde{\QQ}^{(0)}_zdz
\nonumber\\&=\int_{-\infty}^{+\infty}\big(\widetilde{\QQ}_{zzz}^{(0)}-f'(\widetilde{\QQ}^{(0)})\widetilde{\QQ}^{(0)}_z\big):\widetilde{\QQ}^{(2)} dz=0.\nonumber
\end{align}
Thanks to $\widetilde{\QQ}^{(0)}=s(z)(\nn(t, x)\nn(t,x)-\frac{1}{3}\II)$, we find
\ben
&&\int_{-\infty}^{+\infty}\Delta_x\widetilde{\QQ}^{(0)}:\widetilde{\QQ}^{(0)}_z+\nabla_x\widetilde{\QQ}^{(0)}:\na_x\widetilde{\QQ}^{(0)}_zdz=0.\label{equ:4.43}
\een
Hence by (\ref{equ:q-44}) and (\ref{equ:4.43}), we have
\ben
\int_{-\infty}^{+\infty}\big(\widetilde{\HH}^{(-1)}:\widetilde{\QQ}^{(1)}_z
+\widetilde{\HH}^{(0)}:\widetilde{\QQ}^{(0)}_z\big)dz
=\int_{-\infty}^{+\infty}\big(\widetilde{\QQ}^{(0)}_{z}\Delta\varphi+2\nabla\varphi\cdot\nabla_x\widetilde{\QQ}^{(0)}_{z}\big):\widetilde{\QQ}^{(1)}_z.\label{equ:4.44}
\een

Summing up (\ref{equ:q-41})--(\ref{equ:4.41}) and (\ref{equ:4.44}), we conclude that
\begin{align}
&\big[\mathbf{p}^{(0)}\big]-\bigg[\Big\langle2{\mathbf{D}}^{(0)}+\QQ^{(0)}\cdot\HH^{(0)}-\HH^{(0)}\cdot\QQ^{(0)}
-\mathbf{S}_{{\QQ}^{(0)}}({\mathbf{H}}^{(0)})+\sigma^d(\QQ^{(0)},\QQ^{(0)}),{\nabla\varphi}\otimes{\nabla\varphi}\Big\rangle\bigg]
\nonumber\\&=\int_{-\infty}^{+\infty}\Big(\nabla_x\cdot\big(\widetilde{\QQ}^{(0)}\widetilde{\HH}^{(-1)}-\widetilde{\HH}^{(-1)}\widetilde{\QQ}^{(0)}
-\mathbf{S}_{\widetilde{\QQ}^{(0)}}(\widetilde{\HH}^{(-1)})\big)\Big)\cdot{\nabla\varphi}dz
\nonumber\\&\quad-\int_{-\infty}^{+\infty}\big(\widetilde{\QQ}^{(0)}_{z}\Delta\varphi+2\nabla\varphi\cdot\nabla_x\widetilde{\QQ}^{(0)}_{z}\big):\widetilde{\QQ}^{(1)}_z+
\big(\widetilde{\QQ}^{(1)}_z\Delta\varphi+2\na\varphi\cdot\na_x\widetilde{\QQ}^{(1)}_z\big):\widetilde{\QQ}^{(0)}_zdz
\nonumber\\&\quad-\int_{-\infty}^{+\infty}\big(\widetilde{\HH}^{(-1)}:\nabla_x\widetilde{\QQ}^{(0)}\big)\cdot{\nabla\varphi}dz.\label{equ:q-4.41}
\end{align}
By the calculations in section 4.1, the left hand side of (\ref{equ:q-4.41}) exactly equals to
 \begin{align}
\big[\mathbf{p}^{(0)}\big]-\Big[\big\langle\sigma^L+\sigma^E,\nabla\varphi\otimes\nabla\varphi\big\rangle\Big],
\end{align}
where $\sigma^L_-={\mathbf{D}}^{(0)}$ and $\sigma^E_-=0$.

In order that the velocity is continuous across the sharp interface, we take the constant $\xi=0$(see Remark \ref{rem:velocity}).
Then the evolution equation (\ref{equ:MCF-f}) of the sharp interface is reduced to
\begin{align}
\varphi_t&=\Delta\varphi-\mathbf{v}^{(0)}\cdot\nabla\varphi.\nonumber
\end{align}
Hence, the equation (\ref{equ:Q-exp-1}) can be reduced to
\begin{align}
\widetilde{\QQ}^{(1)}_{zz}-f'(\widetilde{\QQ}^{(0)})\widetilde{\QQ}^{(1)}+2\nabla\varphi\cdot\nabla\widetilde{\QQ}^{(0)}_z=
-\mathbf{\widetilde{\Omega}}^{(-1)}\cdot\widetilde{\QQ}^{(0)}+\widetilde{\QQ}^{(0)}\cdot\mathbf{\widetilde{\Omega}}^{(-1)},\nonumber
\end{align}
and
\beno
&&\widetilde{\HH}^{(-1)}=\widetilde{\QQ}^{(0)}_{z}\Delta\varphi-\mathbf{\widetilde{\Omega}}^{(-1)}\cdot\widetilde{\QQ}^{(0)}
+\widetilde{\QQ}^{(0)}\cdot\mathbf{\widetilde{\Omega}}^{(-1)},\\
&&\widetilde{\sigma}^a_{(-1)}=\widetilde{\QQ}^{(0)}\cdot\widetilde{\QQ}^{(0)}\cdot\mathbf{\widetilde{\Omega}}^{(-1)}
+{\widetilde{\Omega}}^{(-1)}\cdot\widetilde{\QQ}^{(0)}\cdot\widetilde{\QQ}^{(0)}.
\eeno
Then by (\ref{equ:v-exp-new0}), we get
\beno
\widetilde{\mathbf{v}}^{(0)}_{z}=-\widetilde{\sigma}^a_{(-1)}\cdot\na\varphi=\f {s(z)} 2\Big((\nn\cdot\na \varphi)^2\widetilde{\mathbf{v}}^{(0)}_{z}
-(\nn\cdot \widetilde{\mathbf{v}}^{(0)}_{z})\big(\nn-\nn\cdot\na \varphi\na\varphi\big)\Big).
\eeno
This implies that $\widetilde{\mathbf{v}}^{(0)}_{z}=0$. Hence, $\big[\vv^{(0)}\big]=0$ and
\beno
\widetilde{\QQ}^{(1)}_{zz}-f'(\widetilde{\QQ}^{(0)})\widetilde{\QQ}^{(1)}+2\nabla\varphi\cdot\nabla\widetilde{\QQ}^{(0)}_z=0.
\eeno
In a similar derivation to (\ref{equ:n-neumann}), we have
\begin{align}
\nu\cdot\nabla\nn=0\quad \textrm{on}\quad \Gamma(t),\label{equ:n-neumann1}
\end{align}
which in turn implies
\begin{align}
-\widetilde{\QQ}^{(1)}_{zz}+f'(\widetilde{\QQ}^{(0)})\widetilde{\QQ}^{(1)}=0.\nonumber
\end{align}
Thus, the solution $\widetilde{\QQ}^{(1)}$ may take
\ben
\widetilde{\QQ}^{(1)}=\pa_z\widetilde{\QQ}^{(0)}\quad \textrm{or}\quad \widetilde{\QQ}{^{(1)}}=s(z)(\nn\nn^\bot+\nn^\bot\nn)\quad \nn^\bot\in \textbf{V}_\nn,\label{equ:Q1-formula}
\een
or $\widetilde{\QQ}^{(1)}=\pa_z\widetilde{\QQ}^{(0)}+s(z)(\nn\nn^\bot+\nn^\bot\nn).$

Therefore, the jump of the pressure becomes
\begin{align}
&\big[\mathbf{p}^{(0)}\big]-\Big[\big\langle\sigma^L+\sigma^E,\nabla\varphi\otimes\nabla\varphi\big\rangle\Big]\nonumber\\
&=-\int_{-\infty}^{+\infty}\big(\widetilde{\QQ}^{(0)}_{z}\Delta\varphi+2\nabla\varphi\cdot\nabla_x\widetilde{\QQ}^{(0)}_{z}\big):\widetilde{\QQ}^{(1)}_z+
\big(\widetilde{\QQ}^{(1)}_z\Delta\varphi+2\na\varphi\cdot\na_x\widetilde{\QQ}^{(1)}_z\big):\widetilde{\QQ}^{(0)}_zdz
\nonumber\\&\quad-\int_{-\infty}^{+\infty}\big(\widetilde{\HH}^{(-1)}:\nabla_x\widetilde{\QQ}^{(0)}\big)\cdot{\nabla\varphi}dz.\nonumber
\end{align}
By (\ref{equ:n-neumann1}) and (\ref{equ:Q1-formula}), we find
\beno
&&\Big[\big\langle\sigma^E,\nabla\varphi\otimes\nabla\varphi\big\rangle\Big]=0,\\
&&\int_{-\infty}^{+\infty}\nabla\varphi\cdot\nabla_x\widetilde{\QQ}^{(0)}_{z}:\widetilde{\QQ}^{(1)}_z+\na\varphi\cdot\na_x\widetilde{\QQ}^{(1)}_z:\widetilde{\QQ}^{(0)}_zdz=0,\\
&&\int_{-\infty}^{+\infty}\big(\widetilde{\HH}^{(-1)}:\nabla_x\widetilde{\QQ}^{(0)}\big)\cdot{\nabla\varphi}dz=0,\\
&&\int_{-\infty}^{+\infty}\big(\widetilde{\QQ}^{(0)}_{z}:\widetilde{\QQ}^{(1)}_z\big)\Delta\varphi+
\big(\widetilde{\QQ}^{(1)}_z:\widetilde{\QQ}^{(0)}_z\big)\Delta\varphi dz=0.
\eeno
Then the jump condition of the pressure is reduced to
\begin{align}
&\big[\mathbf{p}^{(0)}\big]=\Big[\big\langle\sigma^L,\nabla\varphi\otimes\nabla\varphi\big\rangle\Big],\nonumber
\end{align}
where
\beno
\sigma^L=s_+^2\nn\NN-s_+^2\NN\nn+\DD^{(0)},
\eeno
while the symmetric parts induced by $\mathbf{S}_{\QQ^{\varepsilon}}(\DD^{\varepsilon})$ vanish.

\begin{Remark}\label{rem:velocity}
In the case of the constant $\xi\neq 0$, $\widetilde{\mathbf{v}}^{(0)}_{z}\neq 0$. Otherwise, we could conclude from (\ref{equ:v-exp-new1}) that
\beno
\mathbf{S}_{\widetilde{\QQ}^{(0)}}(\widetilde{\QQ}^{(0)})=0,
\eeno
which is impossible. Therefore, the velocity is in general not continuous across the sharp interface in the case of $\xi\neq 0$ except that $\xi=O(\varepsilon)$.

In addition, the assumption on $\xi=0$ also plays an important role for the global well-posedness of the 2-D Beris-Edwards system proved by
Paicu and Zarnescu \cite{Pai}.

\end{Remark}

\section*{Acknowledgments}

M. Fei is partly supported by NSF of China under Grant 11301005.
W. Wang is supported by China Postdoctoral Science Foundation under Grant 2013M540010 and 2014T70008.
P. Zhang is partly supported by NSF of China under Grant 11421101 and 11421110001.
Z. Zhang is partially supported by NSF of China under Grant 11371037, Program for New Century Excellent Talents in University.

\end{document}